\documentclass[journal]{IEEEtran}
\usepackage{soul}
\usepackage{color}
\usepackage[usenames,dvipsnames]{xcolor}

\usepackage{cite}
\usepackage{url}
\usepackage[hidelinks]{hyperref}

\ifCLASSINFOpdf
\else
\fi
\usepackage{tikz}
\usetikzlibrary{arrows, positioning, shapes.geometric, fit, calc,backgrounds,shapes.arrows}

\usepackage{amsmath}
\usepackage{amssymb}
\usepackage{bm}
\usepackage[short]{optidef}
\ifCLASSOPTIONcompsoc
  \usepackage[caption=false,font=normalsize,labelfont=sf,textfont=sf]{subfig}
\else
  \usepackage[caption=false,font=footnotesize]{subfig}
\fi

\renewcommand{\baselinestretch}{0.95}

\begin{document}
%
\title{Learning Active Constraints to Efficiently Solve Linear Bilevel Problems: Application to the Generator Strategic Bidding Problem}

\author{Eléa~Prat
        and Spyros~Chatzivasileiadis,~\IEEEmembership{Senior Member,~IEEE}
\vspace{-0.3cm}
\thanks{E. Prat and S. Chatzivasileiadis are with the Department of Electrical Engineering, Technical University of Denmark, 2800 Kgs. Lyngby, Denmark e-mail: \{emapr, spchatz\}@elektro.dtu.dk.}
\thanks{This work is supported by the H2020 European Project FLEXGRID, Grant Agreement No. 863876  and by the ERC Starting Grant VeriPhIED, Grant Agreement No. 949899}}

\maketitle

\begin{abstract}
Bilevel programming can be used to formulate many problems in the field of power systems, such as strategic bidding. However, common reformulations of bilevel problems to mixed-integer linear programs make solving such problems hard, which impedes their implementation in real-life. In this paper, we  significantly improve solution speed and tractability by introducing decision trees to learn the active constraints of the lower-level problem, while avoiding to introduce binaries and big-M constants. The application of machine learning reduces the online solving time,  by moving the selection of active constraints to an offline process, and becomes particularly beneficial when the same problem has to be solved multiple times. We apply our approach to the strategic bidding of generators in electricity markets, where generators solve the same problem many times for varying load demand or renewable production. Three methods are developed and applied to the problem of a strategic generator, with a DCOPF in the lower-level. These methods are heuristic and as so, do not provide guarantees of optimality or solution quality. Yet, we show that for networks of varying sizes, the computational burden is significantly reduced, while we also manage to find solutions for strategic bidding problems that were previously intractable.

\end{abstract}

\begin{IEEEkeywords}
Bilevel programming, Stackelberg games, classifier, active set, mixed-integer linear programming (MILP), strategic bidding
\end{IEEEkeywords}

%
\IEEEpeerreviewmaketitle

\section{Introduction}
Bilevel problems were formulated for the first time in 1934 by H.v. Stackelberg~\cite{dempe2006foundations}. Since then, they have been widely used in economics and game theory, in particular to model strategic behaviors. One issue is that these problems are NP-hard to solve~\cite{dempe2006foundations, colson2005bilevel}. Linear bilevel problems can easily be reformulated as one-level problems, but the introduction of binary variables renders these reformulations intractable for large systems~\cite{pozo2017bilevel}. In power systems, bilevel problems can be used to model the behavior of a price-maker in electricity markets, to evaluate investment in production facilities, to model the best transmission network investments~\cite{gabriel2012complementarity}, to evaluate the vulnerability of power systems to deliberate~\cite{arroyo2005terrorist} or unintentional~\cite{arroyo2010vulnerability} outages, and more recently for demand response management by tariff design in a smart grid setup~\cite{kovacs2019demand}. Due to the size of the networks, the tractability of bilevel problems is critical for these applications. 

The approach used here derives from the active-set strategy~\cite{fukushima2002active}. The lower-level problem is replaced with its active constraints in a one-level reformulation, avoiding the use of binary variables and thus obtaining a more tractable version of the bilevel problem. However, there can be multiple possible sets of active constraints to consider, depending on the value of the variables of the upper-level problem. To avoid having to identify the possible active sets at every run of a model, machine learning techniques can be used. This is particularly interesting for power systems applications in which similar calculations have to be carried out very often with only a few parameters changing, especially in problems related to bidding in the electricity markets.

In this paper, we study the problem of a strategic generator optimizing its bids on the day-ahead market, with the market modelled by a DC Optimal Power Flow (DCOPF). This problem has been largely studied. It has been formulated for different types of generators, including virtual power plants~\cite{shafiekhani2019} or a retailer considering demand response~\cite{sharifi2020}. Different set-ups have been considered, such as stochastic ones~\cite{ruiz2009}, multi-period, with non-convex operating constraints~\cite{ye2019} or incomplete information~\cite{li2005}. The problem of strategic bidding in different markets at once has also been studied~\cite{rintamaki2020}. However the size of the test cases remains small due to the complexity of solving bilevel problems.

Machine learning has already been used to learn the active constraints of a DCOPF problem with promising results~\cite{ng2018statistical, misra2018learning, deka2019learning}. In~\cite{chen2020learning}, a similar approach is considered, except that the active constraints are not learned directly. Instead, they are derived from the gradient of the cost with respect to the loads; the gradient is itself the output of a Neural Network classifier. In~\cite{xavier2020learning}, machine learning is used to identify redundant constraints and simplify the security-constrained unit commitment problem. Contrary to these problems, the bilevel problems we consider in this paper have two distinct characteristics. Besides being effectively formulated as MILPs -- and not LPs as most of the problems considered above -- their main challenge is that the decision variables of the upper-level problem shall not intervene with the active set classification process of the lower level problem, although the upper-level decision variables are indeed parameters of the lower-level problem. To the best of our knowledge, this is the first paper that introduces machine learning techniques to boost the runtime and solution quality of problems with such features.

Looking at the literature for approaches related to a more efficient solving of linear bilevel problems, several directions have been explored, such as genetic algorithms~\cite{calvete2008genetic} and evolutionary algorithms~\cite{hecheng2012evolutionary, sinha2013efficient}. In~\cite{pineda2018efficient}, regularization approaches are combined with mixed-integer reformulation, by first finding a local optimal solution to provide initial values of the binary variables, which reduces the computational burden. In~\cite{kleinert2021computing}, the problem is reformulated as a one-level problem, using the dual of the lower-level problem, and decomposition is used for solving it.
In these regards, and as shown by the results obtained in this paper, the application of machine learning techniques is a most promising new approach.

This paper has the following contributions: 
\begin{itemize}
  \item We introduce three methods, with some of them being highly parallelizable, that boost the runtime and solution quality of bilevel problems.
  \item Using Decision Trees, we move the selection of the active set of the lower-level problem to an offline process. We eliminate the binary variables, and solve instead a single or a small number of LPs. 
  \item We apply our methods to the problem of the strategic bidding of a generator in the electricity market, and demonstrate their performance to power systems of varying size and complexity, up to 2'869 buses.
  \item We compare our methods with the most promising existing techniques for solving bilevel programs, such as the Big-M method and the penalty alternating direction method (PADM) introduced in~\cite{kleinert2021computing}. We show that our methods are 6-24 times faster while achieving good solution quality, comparable with existing methods, even though our methods do not provide guarantees in these regards, similar to other existing heuristic to solve these non-convex problems. More importantly, our methods are shown to retrieve good solutions to problems that existing methods find intractable.
\end{itemize}

The rest of this paper is organized as follows: Section~\ref{sec:Blvl} introduces bilevel problems as well as the example considered in the rest of the paper with its reformulation as a mixed-integer linear problem (MILP). Section~\ref{sec:meth} describes the methods proposed. The application of the methods to test systems is given in Section~\ref{sec:case}, and Section~\ref{sec:ccl} concludes the paper.

\section{Solving Bilevel Problems}
\label{sec:Blvl}
\subsection{Formulation of KKTs and linearization}
Bilevel problems are optimization problems in which constraints are in part defined by another optimization problem. One common example in the fields of economics and game theory is Stackelberg games, in which one player, the leader, anticipates the decision of the other agents, or followers, and decides on its strategy accordingly. They can be formulated as:
\begin{subequations}
\allowdisplaybreaks
\begin{alignat}{2}
\underset{x,y, \lambda, \mu}{\text{min}} \quad & F(x,y, \lambda, \mu) \\
\text{s.t.} \quad & H(x,y, \lambda, \mu) = 0 \\
& G(x,y, \lambda, \mu) \leq  0\\
& \label{lowerobj} \underset{y}{\text{min}} \quad  && f(x,y) \\
& \label{lowereq} \text{s.t.} \quad && h(x,y) = 0 \quad (\lambda)\\
&&&  \label{lowerineq} g(x,y) \leq  0 \quad (\mu)
\end{alignat}
\end{subequations}
where $x \in \mathbb{R}^n $ and $y \in \mathbb{R}^m $. Equations~\eqref{lowerobj}~to~\eqref{lowerineq} describe the embedded problem, referred to as lower-level or follower problem. The lower-level objective function is $f(x,y)$, and $h(x,y)$ and $g(x,y)$ are the lower-level constraints. The dual variables associated with these equality and inequality constraints respectively are $\lambda$ and $\mu$.
$F(x,y)$ is the objective function of the global problem, called upper-level or leader problem. $H(x,y)$ and $G(x,y)$ are the upper-level constraints. The variables collected in $x$ are decision variables to the upper-level problem and parameters in the lower-level problem. On the other hand, $y$ stands for the decision variables of the lower-level problem.

The problem is non-linear and intractable, due to Equation~\eqref{lowerobj}. In order to be solved, it can be reformulated as a one-level problem. The most common approach, under the condition that the lower-level problem is convex and regular, is to replace it with its Karush--Kuhn--Tucker conditions (KKTs). The following problem is obtained:
\begin{subequations}
\allowdisplaybreaks
\begin{alignat}{2}
\underset{x,y,\lambda,\mu}{\text{min}} \quad & F(x,y, \lambda, \mu) \\
\text{s.t.} \quad & H(x,y, \lambda, \mu) = 0 \\
& G(x,y, \lambda, \mu) \leq  0\\
    \begin{split}
        \label{lderiv} \nabla_y \mathcal{L}(x,y,\lambda,\mu) = \nabla_y f(x,y) + \lambda ^T \nabla_y h(x,y) \\ + \mu ^T \nabla_y g(x,y) = 0
    \end{split} 
    \\
& \label{lowereq2} h(x,y) = 0 \\
& \label{cc2} 0 \leq \mu 	\perp - g(x,y) \geq  0
\end{alignat}
\end{subequations}
where $\nabla_y \mathcal{L}(x,y,\lambda,\mu)$ represents the Lagrangian derivatives with regard to the components of the vector $y$.

The complementarity constraints associated with the inequality constraints of the lower-level problem (Equation~\eqref{cc2}) are non-linear but several techniques exist to linearize them, such as the Fortuny-Amat--McCarl linearization, which will be detailed here. It introduces binary variables, thus transforming the problem into a Mixed Integer Problem (MIP). Equation~\eqref{cc2} can be replaced by:
    \begin{subequations}
    \allowdisplaybreaks
    \begin{alignat}{2}
    \label{cc3} 0 & \leq - g( x,y) \leq  Mu \\
    \label{cc4} 0 & \leq \mu \leq  M(1-u)
    \end{alignat}
    \end{subequations}
    where $M$ is a large enough constant and $u$ is a binary variable. The selection of $M$ is an important issue, which will be discussed in the case studies, in Section~\ref{sec:case}.

If the lower-level problem is linear and the objective function and constraints of the upper-level problem are linear too, the reformulated problem is a MILP and can generally be solved. However, when the number of inequality constraints in the lower-level problem is large, the number of binaries introduced by this reformulation will be high and the problem often becomes intractable. The focus of this paper is, thus, to eliminate these binary variables in order to significantly decrease the solving time and to enable solving problems that were intractable before.

\subsection{Strategic Generator as a MILP}
\label{sec:MILP}
In the rest of this paper, the methods proposed will be applied to one particular instance of bilevel problem which formulates the decision-making of a strategic producer aiming at determining its bids on the day-ahead market, in order to maximize its profit. The lower-level problem is the day-ahead market clearing formulated as a DCOPF:

\begin{small}
\begin{subequations}
\allowdisplaybreaks
\begin{alignat}{2}
\label{strat_lobj} \underset{P^\text{g},\theta}{\text{min}} \quad & c^\text{S} P_{i=1}^\text{g} + \sum_{i\neq 1} c_i P_i^\text{g} \\
\label{strat_bal} \quad \text{s.t.} \quad & P_i^\text{g}-P_i^\text{d}- \sum_{l, i\in l} B_l \Delta\theta_{l} = 0, \quad  \forall i  \quad (\alpha_i)\\
& \label{strat_gen} P_i^\text{g,min} \leq P_i^\text{g} \leq P_i^\text{g,max} , \quad \forall i \quad (\phi_i^\text{min},\phi_i^\text{max})\\
& \label{strat_line} -f_{l}^\text{max} \leq  B_l \Delta\theta_{l} \leq f_{l}^\text{max} , \quad \forall l \quad (\rho_l^\text{min},\rho_l^\text{max})\\
\label{strat_ref} & \theta_{i=\text{ref}}=0 \quad (\gamma)
\end{alignat}
\end{subequations}
\end{small}
where $i \in I$, represents the bus of the system studied and $l \in L$, the lines connecting the bus of this system. The decision variables of the DCOPF are the power output of all the generators in the system $P_i^\text{g}$ $(i \in I)$ and the voltage angles at the bus $\theta_{i}$ $(i \in I)$. $\Delta \theta_{l} $ is a notation to represent the voltage angle difference between the origin bus ($i = \text{from}_l$) and the destination bus ($i = \text{to}_l$) of line $l$, namely, $\Delta \theta_{l}= \theta_{i = \text{from}_l} - \theta_{i = \text{to}_l}$. Without loss of generality, we consider that there is only one generator per bus and that the strategic generator is placed in bus $i=1$. The slack bus is identified as $i = \text{ref}$. It might or might not be the bus where the strategic generator is located. The price bid of the strategic generator $c^\text{S}$ is an upper-level decision and a parameter to the DCOPF.

The objective is to minimize the total cost of the system. The actual production cost of generator $i$ is given by $c_{i}$. Apart from the strategic generator which bids at the cost $c^\text{S}$, all the other generators are assumed to bid their true cost. Equation~\eqref{strat_bal} is the power balance at bus $i$, $B_l$ being the susceptance of line $l$ and $P_i^\text{d}$ the demand at bus $i$. Equation~\eqref{strat_gen} gives the operating limits for the generator at bus $i$ in terms of minimum $P_i^\text{g,min}$ and maximum $P_i^\text{g,max}$. Equation~\eqref{strat_line} gives the limits of the power flow in line $l$, bounded by the line constraint $f_{l}^\text{max}$. Finally, the angle of the slack bus $\theta_{i=\text{ref}}$ is set to 0. The dual variables of Equations~\eqref{strat_bal}~to~\eqref{strat_ref} are given in parentheses next to each equation, and will be used to formulate the KKTs of the lower-level problem.

The bilevel problem for the strategic generator can be formulated as:

\begin{small}
\begin{subequations}
\allowdisplaybreaks
\begin{alignat}{2}
\label{strat_uobj} \underset{c^\text{S}, P^\text{g},\theta, \alpha_{i=1}}{\text{min}} \quad & c_{i=1} P_{i=1}^\text{g}  - \alpha_{i=1} P_{i=1}^\text{g}  \\
\label{strat_cost} \quad \text{s.t.} \quad & c_{i=1} \leq  c^\text{S}  \leq c^\text{S,max} \\
\label{strat_ll} & \eqref{strat_lobj} - \eqref{strat_ref}
\end{alignat}
\end{subequations}
\end{small}
The upper-level problem objective~\eqref{strat_uobj} is the maximization of profit for the strategic generator, as the difference between its operating cost $c_{i=1}$ and the price received, given by the dual variable of the power balance in bus 1, $\alpha_{i=1}$. The profit is here expressed in the minimization form (standard form). Equation~\eqref{strat_cost} belongs to the upper-level problem and sets limits to the cost of the strategic generator to ensure that the problem is bounded. The minimum is equal to the real cost of production $c_{i=1}$ and the maximum is $c^\text{S,max}$.

The objective function in Equation~\eqref{strat_uobj} is not linear because of the term $\alpha_{i=1} P_{i=1}^\text{g}$. However it can be linearized, as shown in~\cite{gabriel2012complementarity}. This together with the KKTs and Fortuny-Amat--McCarl linearization gives the following MILP:

\begin{footnotesize}
\begin{subequations}
\label{eq:milp}
\allowdisplaybreaks
\begin{alignat}{2}
\nonumber \min_{\substack{c^\text{S}, P^\text{g},\theta, \alpha, \\ \rho,\phi,\gamma, u , y}} \quad & \sum_{i\neq 1} (c_i  P_i^\text{g} + \phi_i^\text{max} P_i^\text{g,max} - \phi_i^\text{min} P_i^\text{g,min} - \alpha_i P_i^\text{d}) \\
\label{milp_obj}&  + c_{i=1} P_{i=1}^\text{g} + \sum_{l} f_{l}^\text{max} (\rho_l^\text{min}+\rho_l^\text{max})\\
\label{milp_cost}\quad \text{s.t.} \quad & c_{i=1} \leq  c^\text{S}  \leq c^\text{S, max} \\
\label{milp_bal}& P_i^\text{g}-P_i^\text{d}- \sum_{l, i\in l} B_l \Delta\theta_{l} = 0, \quad  \forall i\\
\label{milp_ref}& \theta_{i=\text{ref}}=0 \\
\label{milp_ls}& c^\text{S} - \alpha_{i=1} - \phi_{i=1}^\text{min} + \phi_{i=1}^\text{max} = 0 \\
\label{milp_lg}& c_i - \alpha_i - \phi_i^\text{min} + \phi_i^\text{max} = 0 \quad \forall i \neq 1\\
\nonumber & \sum_{l, i = \text{from}_l} B_l (\alpha_i - \alpha_{i=\text{to}_l} - \rho_l^\text{min} + \rho_l^\text{max}) \\
\label{milp_lt} & + \sum_{l, i = \text{to}_l} B_l (\alpha_i - \alpha_{i=\text{from}_l} + \rho_l^\text{min} - \rho_l^\text{max}) = 0, \quad \forall i \neq \text{ref}\\
\nonumber & \sum_{l, \text{ref} = \text{from}_l} B_l (\alpha_{i=\text{ref}} - \alpha_{i=\text{to}_l} - \rho_l^\text{min} + \rho_l^\text{max}) \\
\label{milp_lref} & + \sum_{l, \text{ref} = \text{to}_l} B_l (\alpha_{i=\text{ref}} - \alpha_{i=\text{to}_l} + \rho_l^\text{min} - \rho_l^\text{max}) + \gamma = 0 \\
\label{milp_Fgmin}& 0 \leq P_i^\text{g} - P_i^\text{g,min} \leq M u_i^\text{min} , \quad \forall i \\
\label{milp_dgmin}& 0 \leq \phi_i^\text{min} \leq M (1-u_i^\text{min}) , \quad \forall i \\
\label{milp_Fgmax}& 0 \leq P_i^\text{g,max} - P_i^\text{g} \leq M u_i^\text{max} , \quad \forall i \\
\label{milp_dgmax}& 0 \leq \phi_i^\text{max} \leq M (1-u_i^\text{max}) , \quad \forall i \\
\label{milp_Flmin}& 0 \leq f_{l}^\text{max} + B_l \Delta\theta_{l} \leq M y_l^\text{min} , \quad \forall l \\
\label{milp_dlmin}& 0 \leq \rho_l^\text{min} \leq M (1-y_l^\text{min}) , \quad \forall l \\
\label{milp_Flmax}& 0 \leq f_{l}^\text{max} - B_l \Delta\theta_{l} \leq M y_l^\text{max} , \quad \forall l \\
\label{milp_dlmax}& 0 \leq \rho_l^\text{max} \leq M (1-y_l^\text{max}) , \quad \forall l \\
\label{milp_bin}& u_i^\text{min}, u_i^\text{max}, y_l^\text{min}, y_l^\text{max} \in \{0,1\}
\end{alignat}
\end{subequations}
\end{footnotesize}
where  $u_i^\text{min}$, $u_i^\text{max}$, $y_l^\text{min}$ and $y_l^\text{max}$ ($i \in I$, $l \in L$) are the binary variables introduced by Fortuny-Amat McCarl linearization, and $M$ is a large enough constant. Equations~\eqref{milp_bal}~and~\eqref{milp_ref} are the equalities of the lower-level problem. Equations~\eqref{milp_ls}~to~\eqref{milp_lref} are obtained by setting to zero the derivatives of the Lagragian of the lower-level problem with regard to all the variables. Equations~\eqref{milp_Fgmin}~to~\eqref{milp_bin} are the linearized complementarity constraints. They contain and replace the inequality constraints in Equations~\eqref{strat_gen}~and~\eqref{strat_line}.

Solving Equations~\eqref{milp_obj}~to~\eqref{milp_bin} directly will be used as a baseline for the case studies in Section~\ref{sec:case}. 

\section{Methods}
\label{sec:meth}
The reformulation proposed here is based on the model given in~\ref{sec:MILP}, and aims at only keeping the constraints that are active at the optimal point. In the case of a linear problem, those are sufficient to describe the system at optimality. We have established three methods to achieve this. The process for each of these methods is illustrated in Figure~\ref{fig:flowcharts}. They follow the same general structure. First, as part of an offline process, a database is built, mapping the variables of the lower-level problem to the corresponding active constraints. This is described in Section~\ref{sec:DBbuild}. This database is used to train a decision tree (DT), as explained in Section~\ref{sec:DTbuild}. For a given value of the parameters, this DT allows to retrieve sets of active constraints, in order to build and solve a reduced bilevel problem. This process is detailed in Section~\ref{sec:RedBlvl}. A  summary of the learning process for the three methods is given in Table~\ref{tab:details}. Finally, a preliminary discussion is available in Section~\ref{sec:Disc}.

\begin{table}[b]
\caption{Comparison of the three methods introduced, in terms of database, DT and reduced bilevel problem}
\resizebox{0.5\textwidth}{!}{
\begin{tabular}{lllllll}
\hline
       & \multicolumn{2}{l}{Database} & \multicolumn{3}{l}{Final DT(s)} & Online       \\
Method & Features       & Target      & DT(s)   & Features   & Target   & LP(s) \\ \hline
VarLower &
  Load and $c^\text{S}$ &
  \begin{tabular}[c]{@{}l@{}}Set of active\\ constraints\end{tabular} &
  $n_\text{V}$ &
  Load &
  \begin{tabular}[c]{@{}l@{}}Set of active\\ constraints\end{tabular} &
  $n_\text{V}$ \\
AllSets &
  Load &
  \begin{tabular}[c]{@{}l@{}}Set of sets of\\ active constraints\end{tabular} &
  1 &
  Load &
  \begin{tabular}[c]{@{}l@{}}Set of $n_\text{A}$ sets of\\ active constraints\end{tabular} &
  $n_\text{A}$ \\
BestSet &
  Load &
  \begin{tabular}[c]{@{}l@{}}Set of active\\ constraints\end{tabular} &
  1 &
  Load &
  \begin{tabular}[c]{@{}l@{}}Set of active\\ constraints\end{tabular} &
  1 \\ \hline
\end{tabular}}
\label{tab:details}
\end{table}


\tikzstyle{io} = [trapezium, trapezium left angle=80, trapezium right angle=100, minimum width=1cm, minimum height=0.25cm, text centered, draw=black, fill=white]
\tikzstyle{ioGroup} = [trapezium, trapezium left angle=80, trapezium right angle=100, text centered, draw=black, thick]
\tikzstyle{roundGroup} = [rectangle, rounded corners = 10, minimum width=100, minimum height=30,text centered, draw=black, thick]
\tikzstyle{process} = [rectangle, minimum width=1cm, minimum height=1cm, text centered, draw=black, fill=white]
\tikzstyle{processGroup} = [rectangle, minimum width=1cm, minimum height=1cm, text centered, draw=black, thick]
\tikzstyle{database} = [cylinder, fill=white, shape border rotate=90,draw,minimum height=1.5cm, minimum width=2cm,shape aspect=.25,]
\tikzstyle{arrow} = [->,>=stealth]
        
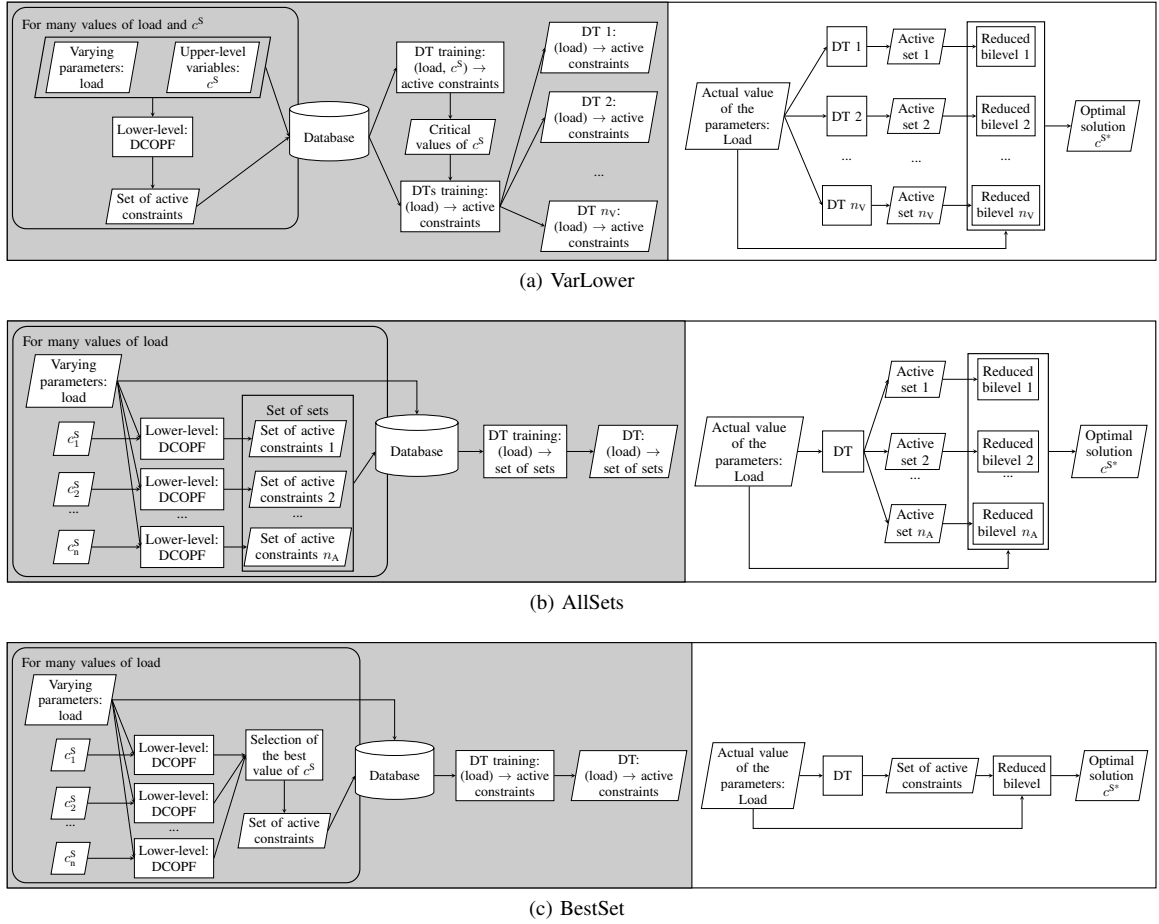
\begin{figure*}[hbt!]
    \centering
    \subfloat[VarLower]{
         \label{fig:meth2}
         \centering
            
        \resizebox{0.85\textwidth}{!}{
        \begin{tikzpicture}[node distance=50]
        
        \node [align=center]  (VarParam) [io] {Varying \\ parameters: \\ load};
        
        \node [align=center] (UpperVar) [io, right of = VarParam,xshift=1.3cm] {
        Upper-level \\ variables: \\ $c^\text{S}$};

        \node [ioGroup,fit=(VarParam)(UpperVar)] [label={[xshift=-1cm, name=l1] For many values of load and $c^\text{S}$}](Input) {};
        
        \node [align=center] (DCOPF) [process, below of = Input] {
        Lower-level: \\ DCOPF } ;
        
        \node [align=center]  (AS) [io, below of = DCOPF] {Set of active \\ constraints};
        
        \node [roundGroup,fit={(l1)(AS)($(Input.east)+(0.7,0)$)}] (DBGen) {};
        
        \draw [arrow] (Input) -- (DCOPF);
        \draw [arrow] (DCOPF) -- (AS);
        
        \node [align=center] (DB) [database, right of = DCOPF,xshift=2.7cm] {Database};
        \draw [arrow] (Input.east) -- (DB.west);
        \draw [arrow] (AS.east) -- (DB.west);
        
        \node [align=center] (critical) [io, right of = DB,xshift=1.3cm] {Critical \\ values of $c^\text{S}$ };
        
        \node [align=center] (DT) [process, above of = critical] {DT training: \\ (load, $c^\text{S}$) $\rightarrow$ \\ active constraints};
        
        \node [align=center] (DTd) [process, below of = critical] {DTs training: \\ (load) $\rightarrow$ active \\ constraints};
        
        \draw [arrow] (DB.east) -- (DT.west);
        \draw [arrow] (DB.east) -- (DTd.west);
        \draw [arrow] (DT) -- (critical);
        \draw [arrow] (critical) -- (DTd);
        
        \node [align=center] (DTd1) [io, right of = DT, xshift=2cm, yshift=0.5cm] {DT 1: \\ (load) $\rightarrow$ active \\ constraints};
        
        \node [align=center] (DTd2) [io, below of = DTd1] {DT 2: \\ (load) $\rightarrow$ active \\ constraints};
        
        \node [align=center][label={[yshift=0.5cm] ...}] (DTdn) [io, below of = DTd2, yshift = -1cm] {DT $n_\text{V}$: \\ (load) $\rightarrow$ active \\ constraints};
        
        \draw [arrow] (DTd.east) -- (DTd1.west);
        \draw [arrow] (DTd.east) -- (DTd2.west);
        \draw [arrow] (DTd.east) -- (DTdn.west);
        
        \node [align=center] (online) [io, right of = critical, xshift = 5.5cm, yshift = 15] {Actual value \\ of the \\ parameters:\\ Load};
        
        \node [align=center] (DT2) [process, right of = online, xshift = 1cm] {DT 2};
        
        \node [align=center] (DT1) [process, above of = DT2, xshift = 0cm] {DT 1};
        
        \node [align=center] [label={[yshift=0.5cm] ...}] (DTn) [process, below of = DT2, yshift = -0.5cm] {DT $n_\text{V}$};
        
        \draw [arrow] (online.east) -- (DT1.west);
        \draw [arrow] (online.east) -- (DT2.west);
        \draw [arrow] (online.east) -- (DTn.west);
        
        \node [align=center] (set1) [io, right of = DT1] {Active \\ set 1};
        
        \node [align=center] (set2) [io, right of = DT2] {Active \\ set 2};
        
        \node [align=center][label={[yshift=0.5cm] ...}] (setn) [io, right of = DTn] {Active \\ set $n_\text{V}$};
        
        \draw [arrow] (DT1.east) -- (set1.west);
        \draw [arrow] (DT2.east) -- (set2.west);
        \draw [arrow] (DTn.east) -- (setn.west);
        
        \node [align=center] (bv1) [process, right of = set1, xshift=0.5cm] {Reduced \\ bilevel 1};
        
        \node [align=center] (bv2) [process, right of = set2, xshift=0.5cm] {Reduced \\ bilevel 2};
        
        \node [align=center][label={[yshift=0.5cm] ...}] (bvn) [process, right of = setn, xshift=0.5cm] {Reduced \\ bilevel $n_\text{V}$};
        
        \draw [arrow] (set1.east) -- (bv1.west);
        \draw [arrow] (set2.east) -- (bv2.west);
        \draw [arrow] (setn.east) -- (bvn.west);
        
        \node [processGroup,fit={(bv1)(bvn)}] (bv) {};
        
        \draw [arrow] (online.south)--++ (0,-2.5) -| (bv.south);
        
        \node [align=center] (sol) [io, right of = bv, xshift=20] {Optimal \\ solution \\ $c^\text{S*}$};
        
        \draw [arrow] (bv.east) -- (sol.west);
        
        \begin{scope}[on background layer]
        \node [rectangle, minimum width=1cm, minimum height=1cm, text centered, draw=black, thick,fit={(DBGen)($(sol.east)+(0.4,0)$)($(DTdn.south)+(0,-0.1)$)}] (All) {};
        \node [rectangle, minimum width=1cm, minimum height=1cm, text centered, draw=black,fill=black!20, thick,fit={(DBGen)($(online.west)+(-0.7,0)$)($(DTdn.south)+(0,-0.1)$)}] (offline) {};
        \end{scope}
        
        \end{tikzpicture}
        }}

    \subfloat[AllSets]{
        \centering
        \label{fig:meth3}
        
        \resizebox{0.85\textwidth}{!}{
        \begin{tikzpicture}[node distance=50]
        
        \node [align=center]  (VarParam) [label={[xshift=15, name=l1] For many values of load}][io] {Varying \\ parameters: \\ load};
        
        \node [align=center] (cs1) [io, below of = VarParam, yshift=10] {$c^\text{S}_1$};
        \node [align=center] (cs2) [io, below of = cs1, yshift=15] {$c^\text{S}_2$};
        \node [align=center][label={[yshift=8] ...}] (csn) [io, below of = cs2, yshift=10] {$c^\text{S}_\text{n}$};
        
        \node [align=center] (DCOPF1) [process, right of = cs1, xshift=25] {Lower-level: \\ DCOPF } ;
        \node [align=center] (DCOPF2) [process, right of = cs2, xshift=25] {Lower-level: \\ DCOPF } ;
        \node [align=center] (DCOPFn) [label={[yshift=2] ...}][process, right of = csn, xshift=25] {Lower-level: \\ DCOPF } ;
        
        \node [align=center]  (AS1) [label={[xshift=0, name=l2] Set of sets}][io, right of = DCOPF1, xshift=30] {Set of active \\ constraints 1};
        \node [align=center]  (AS2) [io, right of = DCOPF2, xshift=30] {Set of active \\ constraints 2};
        \node [align=center]  (ASn) [label={[yshift=4] ...}][io, right of = DCOPFn, xshift=30] {Set of active \\ constraints $n_\text{A}$};
        
        \node [processGroup,fit=(l2)(ASn)] (SAS) {};
        
        \draw [arrow] (VarParam.east) -- (DCOPF1.west);
        \draw [arrow] (VarParam.east) -- (DCOPF2.west);
        \draw [arrow] (VarParam.east) -- (DCOPFn.west);
        \draw [arrow] (cs1.east) -- (DCOPF1.west);
        \draw [arrow] (cs2.east) -- (DCOPF2.west);
        \draw [arrow] (csn.east) -- (DCOPFn.west);
        \draw [arrow] (DCOPF1.east) -- (AS1.west);
        \draw [arrow] (DCOPF2.east) -- (AS2.west);
        \draw [arrow] (DCOPFn.east) -- (ASn.west);
        
        \node [roundGroup,fit={(l1)(SAS)($(SAS.east)+(0.7,0)$)}] (DBGen) {};
        
        \node [align=center] (DB) [database, right of = DBGen,xshift=100] {Database};
        \draw [arrow] (VarParam.east) -| (DB.north);
        \draw [arrow] (SAS.east) -- (DB.west);
        
        \node [align=center] (DT) [process, right of = DB,xshift=25] {DT training: \\ (load) $\rightarrow$ \\ set of sets};
        
        \node [align=center] (DTd) [io, right of = DT,xshift=25] {DT: \\ (load) $\rightarrow$ \\ set of sets};
        
        \draw [arrow] (DB.east) -- (DT.west);
        \draw [arrow] (DT) -- (DTd);
        
        \node [align=center] (online) [io, right of = DTd, xshift = 30, yshift = 0] {Actual value \\ of the \\ parameters:\\ Load};
        
        \node [align=center] (DTo) [process, right of = online, xshift=15] {DT};
        \draw [arrow] (online) -- (DTo);
        
        \node [align=center] (set2) [io, right of = DTo] {Active \\ set 2};
        \node [align=center] (set1) [io, above of = set2] {Active \\ set 1};
        \node [align=center][label={[yshift=0.5cm] ...}] (setn) [io, below of = set2] {Active \\ set $n_\text{A}$};
        
        \draw [arrow] (DTo.east) -- (set1.west);
        \draw [arrow] (DTo.east) -- (set2.west);
        \draw [arrow] (DTo.east) -- (setn.west);
        
        \node [align=center] (bv1) [process, right of = set1, xshift=0.5cm] {Reduced \\ bilevel 1};
        
        \node [align=center] (bv2) [process, right of = set2, xshift=0.5cm] {Reduced \\ bilevel 2};
        
        \node [align=center][label={[yshift=0.5cm] ...}] (bvn) [process, right of = setn, xshift=0.5cm] {Reduced \\ bilevel $n_\text{A}$};
        
        \draw [arrow] (set1.east) -- (bv1.west);
        \draw [arrow] (set2.east) -- (bv2.west);
        \draw [arrow] (setn.east) -- (bvn.west);
        
        \node [processGroup,fit={(bv1)(bvn)}] (bv) {};
        
        \draw [arrow] (online.south)--++ (0,-2) -| (bv.south);
        
        \node [align=center] (sol) [io, right of = bv, xshift=20] {Optimal \\ solution \\ $c^\text{S*}$};
        
        \draw [arrow] (bv.east) -- (sol.west);
        
        \begin{scope}[on background layer]
        \node [rectangle, minimum width=1cm, minimum height=1cm, text centered, draw=black, thick,fit={(DBGen)($(sol.east)+(0.2,0)$)}] (All) {};
        \node [rectangle, minimum width=1cm, minimum height=1cm, text centered, draw=black,fill=black!20, thick,fit={(DBGen)($(online.west)+(-0.5,0)$)}] (offline) {};
        \end{scope}
        
        \end{tikzpicture}
        }}

    \subfloat[BestSet]{
        \centering
        \label{fig:meth4}
        
        \resizebox{0.85\textwidth}{!}{
        \begin{tikzpicture}[node distance=50]
        
        \node [align=center]  (VarParam) [label={[xshift=15, name=l1] For many values of load}][io] {Varying \\ parameters: \\ load};
        
        \node [align=center] (cs1) [io, below of = VarParam, yshift=10] {$c^\text{S}_1$};
        \node [align=center] (cs2) [io, below of = cs1, yshift=15] {$c^\text{S}_2$};
        \node [align=center][label={[yshift=8] ...}] (csn) [io, below of = cs2, yshift=10] {$c^\text{S}_\text{n}$};
        
        \node [align=center] (DCOPF1) [process, right of = cs1, xshift=25] {Lower-level: \\ DCOPF } ;
        \node [align=center] (DCOPF2) [process, right of = cs2, xshift=25] {Lower-level: \\ DCOPF } ;
        \node [align=center] (DCOPFn) [label={[yshift=2] ...}][process, right of = csn, xshift=25] {Lower-level: \\ DCOPF } ;
        
        \node [align=center]  (best) [process, right of = DCOPF1, xshift=30] {Selection of \\ the best \\ value of $c^\text{S}$};
        
        \node [align=center]  (AS) [io, below of = best, yshift=-5] {Set of active \\ constraints};
        
        \draw [arrow] (VarParam.east) -- (DCOPF1.west);
        \draw [arrow] (VarParam.east) -- (DCOPF2.west);
        \draw [arrow] (VarParam.east) -- (DCOPFn.west);
        \draw [arrow] (cs1.east) -- (DCOPF1.west);
        \draw [arrow] (cs2.east) -- (DCOPF2.west);
        \draw [arrow] (csn.east) -- (DCOPFn.west);
        \draw [arrow] (DCOPF1.east) -- (best.west);
        \draw [arrow] (DCOPF2.east) -- (best.west);
        \draw [arrow] (DCOPFn.east) -- (best.west);
        \draw [arrow] (best) -- (AS);
        
        \node [roundGroup,fit={(l1)(DCOPFn)($(AS.east)+(0.7,0)$)}] (DBGen) {};
        
        \node [align=center] (DB) [database, right of = AS,xshift=30, yshift=40] {Database};
        \draw [arrow] (VarParam.east) -| (DB.north);
        \draw [arrow] (AS.east) -- (DB.west);
        
        \node [align=center] (DT) [process, right of = DB,xshift=30] {DT training: \\ (load) $\rightarrow$ active \\ constraints};
        
        \node [align=center] (DTd) [io, right of = DT,xshift=40] {DT: \\ (load) $\rightarrow$ active \\ constraints};
        
        \draw [arrow] (DB.east) -- (DT.west);
        \draw [arrow] (DT) -- (DTd);
        
        \node [align=center] (online) [io, right of = DTd, xshift = 40] {Actual value \\ of the \\ parameters:\\ Load};
        
        \node [align=center] (DTo) [process, right of = online, xshift=15] {DT};
        \draw [arrow] (online) -- (DTo);
        
        \node [align=center] (set) [io, right of = DTo, xshift=15] {Set of active \\ constraints};
        
        \draw [arrow] (DTo.east) -- (set.west);
        
        \node [align=center] (bv) [process, right of = set, xshift=15] {Reduced \\ bilevel };
        
        \draw [arrow] (set.east) -- (bv.west);
        
        \draw [arrow] (online.south)--++ (0,-0.5) -| (bv.south);
        
        \node [align=center] (sol) [io, right of = bv, xshift=15] {Optimal \\ solution \\ $c^\text{S*}$};
        
        \draw [arrow] (bv.east) -- (sol.west);
        
        \begin{scope}[on background layer]
        \node [rectangle, minimum width=1cm, minimum height=1cm, text centered, draw=black, thick,fit={(DBGen)($(sol.east)+(0.2,0)$)}] (All) {};
        \node [rectangle, minimum width=1cm, minimum height=1cm, text centered, draw=black,fill=black!20, thick,fit={(DBGen)($(online.west)+(-0.5,0)$)}] (offline) {};
        \end{scope}
        
        \end{tikzpicture}
        }}
    \caption{Description of the three methods introduced in this paper. The grayed area on the left contains the offline part of the method, while the online part is represented in the right part. The step ``DT training" also includes testing of the decision tree.}
    \label{fig:flowcharts}
\end{figure*}


\subsection{Database Generation and Learning}
\label{sec:learning}
The general idea is to reduce the lower-level problem to its active constraints, which will eliminate the binaries introduced by the linearization of the complementarity constraints. The optimal solution will be different for different values of the input parameters (such as the loads in the case of the strategic generator problem) and so will the active constraints. As a consequence, the identification of these active constraints must be carried out for each new value of the input parameters. This can be tedious as the variables of the upper-level problem are also parameters to the lower-level problem. So, in order to consider all possible reductions of the lower-level problem, multiple setups would have to be tested, even when the parameters are known. In the context of power systems, this is particularly critical as decisions have to be made very often and the parameters vary and are uncertain, especially demand and renewable energy generation. To avoid a long decision process, the idea is to move the selection of the active constraints to an offline process, using machine learning classification techniques. In this paper, DTs are used to perform this classification. 
Our approach can also be adapted to other machine learning classification approaches, such as Neural Networks or Random Forests (see e.g.~\cite{murzakhanov2021} about how we can use an exact transformation of Neural Networks to a MILP for power system problems). Here, we focus on DTs because they are easier to interpret~\cite{carrizosa2011}. For this reason, they have a larger acceptance in the industry, which we think is crucial for the deployment of those methods in actual practice.
The challenge is to remove the decision variables of the upper-level problem from the classifier, unless the DT is included in the optimization problem as in~\cite{halilbasic2018scopf}. This option has been considered but discarded as it would introduce unnecessary binary variables and additional constraints, while preventing the removal of constraints. The three methods detailed in the following offer three different ways to exclude the decision variables of the upper-level problem from the classification. Solving the MILP reformulation for each possible active set of lower-level problem could be a way of proceeding. But while this works well for a small number of active sets, it becomes highly inefficient when there are too many. On the other hand, for a system that has few active sets in the lower-level, it will be more suitable to directly solve the MILP for each of the sets than to use the DT approach.

\subsubsection{Identification of the Active Constraints}
The focus here is on the lower-level optimization problem. In optimization problems, if the feasible set is convex and variables are continuous, the optimal solution lies at the boundaries of the feasible space. In particular, this is the case for a linear problem (LP), as considered here. The corresponding constraints are binding, and the other constraints are inactive. The active set that corresponds to the optimal solution regroups the constraints that are satisfied with equality at the optimal point. If the problem is reformulated following this set, by replacing the active inequality constraints with equalities and removing the other constraints, it will recover the optimal solution. This reduced problem is easier to solve since some constraints are dropped.

The active constraints of an LP can partially be identified by looking at the value of the dual variables associated with the inequalities of the problem, at the optimal point (equality constraints are always binding). All dual variables that are non-zero indicate an active constraint. However, a dual variable equal to zero might also indicate an active constraint. It is then necessary to check the activation by looking at the value of the constraint at the optimal point.

\subsubsection{Database Building Process}
\label{sec:DBbuild}
In order to train the DT classifier, a database of points has to be built. The way this database is generated is different in each method we propose but the general idea is the same: the DT should take the load at every bus as input and return one or several active sets. We need to ensure that the database generated gives a good representation of the possible active sets. To achieve this, the algorithm \textit{DiscoverMass}, as presented in~\cite{misra2018learning}, is applied as a stopping criterion. The idea is to keep generating points from a given distribution, until a sufficient share of the possible active sets have been recovered; that is, until the probability mass of the discovered sets reaches a chosen threshold. A safety limit to the number of steps of the algorithm is also defined, in case it would not converge fast enough. The interested reader can refer to~\cite{misra2018learning} for details on the algorithms and the theorems and proofs associated, in particular regarding termination.

\paragraph{Method VarLower}
In this method, shown in Figure~\ref{fig:meth2}, the generated point consists of the varying parameters of the lower-level and the variables of the upper-level problem that are parameters to the lower-level problem. In the selected example, that would be the all the bus loads, collected in a load vector, along with the cost bid of the strategic generator. For the randomly generated point, the lower-level problem (DCOPF) is run and the set of active constraints at the optimal point is retrieved. 

\paragraph{Method AllSets}
As illustrated in Figure~\ref{fig:meth3}, the database in this method associates for a given load vector all observed sets of active constraints obtained by varying the cost bid of the strategic generator. For each randomly created load vector, multiple instances of the DCOPF are solved for a range of $c^\text{S}$ values. For each value of $c^\text{S}$, an active set is retrieved. All these active sets are gathered in a set of active sets to be associated with this load in the database.

\paragraph{Method BestSet}
This works similarly to the previous method, as shown in Figure~\ref{fig:meth4}: for a given load vector, the DCOPF is solved for a range of $c^\text{S}$ values. In this case, however, only the active set corresponding to the value of $c^\text{S}$ that returns the best value of the upper-level objective function is kept to be part of the database. The intention here is to keep only the active sets corresponding to optimal points of the bilevel problem.

\subsubsection{Decision Tree Training}
\label{sec:DTbuild}
\begin{figure*}[ht!]
    \centering
    \tikzstyle{set} = [rectangle, minimum width=1.5cm, minimum height=0.7cm, text centered, draw=black, fill=black!20]
    \tikzstyle{cs} = [rectangle, minimum width=1.5cm, minimum height=0.7cm, text centered, draw=red, fill=white,text=red, thick]
    \tikzstyle{d} = [rectangle, minimum width=2cm, minimum height=0.7cm, text centered, draw=black, fill=white]
    
    \resizebox{1\textwidth}{!}{
    \begin{tikzpicture}
      [
        grow = down,
        level distance = 1.5cm
      ]
        
        \node (DT1) [d] {Is $P^\text{d}_\text{1} \leq 7.8$  ?} [sibling distance=6.5cm]
            child { [sibling distance=2.5cm] node [cs] {\textbf{Is $\bm{c^\text{S} \leq 10}$  ?}}
                child {[sibling distance=2cm] node [set] {Set 1}}
                child {[sibling distance=2cm]node [d] {Is $P^\text{d}_\text{1} \leq 5$  ?}
                    child {node [set] {Set 5}}
                    child {node [set] {Set 6}}}
                edge from parent node [above] {yes}
            }
            child { [sibling distance=2.5cm] node [d] {Is $P^\text{d}_\text{2} \leq 15$  ?}
                child { [sibling distance=2cm]node [cs] {\textbf{Is $\bm{c^\text{S} \leq 18}$  ?}}
                    child {node [set] {Set 3}}
                    child {node [set] {Set 4}}
                }
                child {[sibling distance=2cm] node [set] {Set 2}}
                edge from parent node [above] {no}
            };
      
        \node (DT2) [label={[xshift=-1.5cm, yshift=0.2cm, name=l2] For $\text{c}^\text{S} \leq 10$}] [d, right of = DT1, xshift = 8cm, yshift = -1.5cm] {Is $P^\text{d}_\text{1} \leq 7.8$  ?} [sibling distance=3cm]
            child {node [set] {Set 1}}
            child {[sibling distance=2cm] node [d] {Is $P^\text{d}_\text{2} \leq 15$  ?}
                child {node [set] {Set 3}}
                child {node [set] {Set 2}}
            };
            
        \node (DT3) [label={[xshift=-1.5cm,yshift=0.2cm, name=l3] For $10 < \text{c}^\text{S} \leq 18$}] [d, right of = DT2, xshift = 6.6cm] {Is $P^\text{d}_\text{1} \leq 7.8$  ?} [sibling distance=4cm]
            child { [sibling distance=2cm] node [d] {Is $P^\text{d}_\text{1} \leq 5$  ?}
                child {node [set] {Set 5}}
                child {node [set] {Set 6}}
                }
            child { [sibling distance=2cm] node [d] {Is $P^\text{d}_\text{2} \leq 15$  ?}
                child {node [set] {Set 3}}
                child {node [set] {Set 2}}
            };
            
        \node (DT4) [label={[xshift=-1.5cm,yshift=0.2cm, name=l4] For $\text{c}^\text{S} > 18$}] [d, right of = DT3, xshift = 7.2cm] {Is $P^\text{d}_\text{1} \leq 7.8$  ?} [sibling distance=4cm]
            child { [sibling distance=2cm] node [d] {Is $P^\text{d}_\text{1} \leq 5$  ?}
                child {node [set] {Set 5}}
                child {node [set] {Set 6}}
                }
            child { [sibling distance=2cm] node [d] {Is $P^\text{d}_\text{2} \leq 15$  ?}
                child {node [set] {Set 4}}
                child {node [set] {Set 2}}
            };
            
        \node [rectangle, minimum width=1cm, minimum height=1cm, draw=black, thick,fit={($(DT1.west)+(-4.3,0.5)$)($(DT4.east)+(2.8,-3.5)$)}] (All) {};
        \draw (5.9,0.6) -- (5.9,-5.1);
        \draw [dashed] (12.55,0.6) -- (12.55,-5.1);
        \draw [dashed] (20.7,0.6) -- (20.7,-5.1);
        
        \node[draw, single arrow,
          minimum height=15mm, minimum width=9mm,
          single arrow head extend=1.7mm,
          anchor=west, fill=white] at (5.2,-2) {};
        
    \end{tikzpicture}
    }
    \caption{Decision tree building process for the method \textit{VarLower}. A first DT is built including load and $c^\text{S}$, in order to identify the critical values of $c^\text{S}$. This helps building intervals of $c^\text{S}$ to split the database and build new decision trees on the load only, one for each interval identified.}
    \label{fig:DT}
\end{figure*}
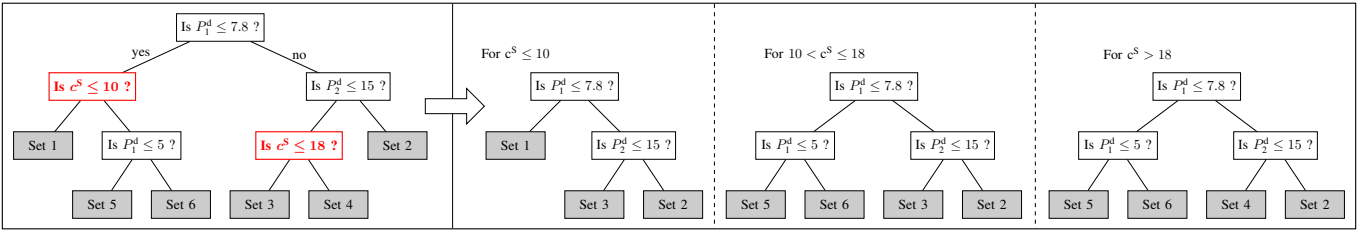

Once the database is created, a DT is trained in order to later predict, for any given value of the parameters of the lower-level problem, the set of active constraints to apply. A decision tree is a classifier that keeps splitting the data according to one of the features of the input until reaching a separation per class. It consists of nodes, which represent decisions made based on a given feature, branches, and leaves, which are the final nodes, in which the class is selected. Figure \ref{fig:DT} shows representations of such decision trees.


\paragraph{Method VarLower}
In this method, a first DT is generated from the database, that is from samples consisting of load vector and $c^\text{S}$ values. However, $c^\text{S}$ should not be an input parameter to the final DT, which will be applied outside and before the bilevel problem is solved. As a consequence, all the nodes in which the feature is $c^\text{S}$ are retrieved and the corresponding critical values of $c^\text{S}$ are extracted, as shown in Figure~\ref{fig:DT}. 
The database is then split following the identified intervals of $c^\text{S}$. Then, for each of these intervals, a DT is built, taking as input the load vector and returning an active set. Online, each of these sub-decision trees will be applied, thus returning as many active sets as there are sub-decision trees. This whole process is illustrated in Figure~\ref{fig:DT}.

\paragraph{Method AllSets}
For the method \textit{AllSets}, the database is already built for the load vector only, so a single DT is built, taking as input the load vector and returning a set of active sets.

\paragraph{Method BestSet}
Similarly to the previous method, one DT is built with the load vector as input, but this time only one active set is returned: the one corresponding to the optimal solution of the bilevel problem.

\subsection{Reduced Bilevel}
\label{sec:RedBlvl}
Once one or more (for method \textit{VarLower}) DTs have been trained offline, they can be applied online to predict the active set corresponding to a given load vector. The lower-level problem is then reduced to keep the active constraints only, replacing them with equalities. This reformulation is applied to the inequalities~\eqref{cc3}~and~\eqref{cc4} in the MILP formulation. For each inequality of the original lower-level problem:
\begin{itemize}
    \item If it is considered as active, it is replaced with equality and the corresponding complementarity constraint is discarded, as well as the constraints on the dual variables, since they are inactive. In the general formulation of the bilevel program, for $g$ active, we would replace Equations~\eqref{cc3}~and~\eqref{cc4} with:
    \begin{equation}
    g(x,y) =  0
    \end{equation}
    \item If it is considered as inactive, the inequality is dropped and the corresponding dual variable is set to be equal to 0. In the general formulation of the bilevel program, for g inactive, we would replace Equations~\eqref{cc3}~and~\eqref{cc4} with:
    \begin{equation}
    \mu =  0
    \end{equation}
\end{itemize}
As a consequence, the binary variables are completely removed from the formulation and the bilevel problem is now an LP. Inequality constraints are also removed. 

Note that we still have the equalities~\eqref{milp_ls}~to~\eqref{milp_lref} associated with the Lagrangian derivatives in the new formulation. In general, for a nondegenerated LP it is possible to simply replace the lower-level problem with its active constraints, without having to formulate any Lagrangian, as in~\cite{chen2020learning}. However, the chosen problem could be degenerated, and moreover, as the duals of the lower-level problem intervene in the objective function~\eqref{milp_obj}, this cannot be applied here.

In the cases where several active sets are returned by the DT(s), several LPs are solved and compared in terms of value of the objective function, to keep the best one only.
\paragraph{Method VarLower}
In the method \textit{VarLower}, there are $\text{n}_\text{V}$ decision trees to be applied, one per interval of $c^\text{S}$ as illustrated in the online part of Figure~\ref{fig:meth2}. As a result, $\text{n}_\text{V}$ sets of active constraints are returned.  The bilevel problem, reformulated as an LP with the corresponding active constraints only, is solved $\text{n}_\text{V}$ times.

\paragraph{Method AllSets}
With \textit{AllSets}, one DT is applied, and the output is a set of $\text{n}_\text{A}$ active sets. $\text{n}_\text{A}$ different versions of the reduced bilevel are solved and the best result is identified. This is shown in the online part of Figure~\ref{fig:meth3}.

\paragraph{Method BestSet}
In the method \textit{BestSet}, one DT is applied, returning one active set. Only one reduced bilevel problem is solved in this case.

\subsection{Short Discussion}
\label{sec:Disc}

The methods introduced are based on the constraints of the problem at hand. In case of network reconfiguration, the learning process has to be applied again. An evolution of these methods would introduce a more flexible setup, similarly to the Graph Neural Network approach followed in~\cite{Donon2019}.

In bilevel problems, it is possible that multiple values of the upper-level variable could be optimal. In this case, one must decide between following an optimistic or a pessimistic approach~\cite{colson2005bilevel}. Here, we consider an optimistic setup: the cost chosen by the generator in the set of possible optimal values does not have an impact on the outcome. However, it could be the case that the strategic generator wants to reduce the risk by bidding the smallest cost in that interval, to ensure being dispatched. In this case, the pessimistic modeling approach would be applied. As the pessimistic approach is not compatible with the KKTs reformulation, our methods could not be used directly. Further research would be necessary on this topic.

In general, the methods introduced can easily be applied to other linear bilevel problems with a unique upper-level decision variable, in an optimistic setup.

\section{Case Studies}
\label{sec:case}
In this section, the three methods will be applied to 5 different test cases, and will be compared with solving the bilevel problem with the KKTs and big-M reformulation as well as with the penalty alternating direction method (PADM) introduced in~\cite{kleinert2021computing}. The systems are based on Matpower cases~\cite{Matpower}. The details of the cases are given in Table~\ref{tab:cases}. The exact data files used are available online~\cite{Zenodo}. All simulations were carried out in Python using Gurobi to solve the optimization problems and Python library scikit-learn to build and apply the decision trees. The code accompanying this paper is also available online~\cite{GitHub}.

\begin{table}[b]
\caption{Characteristics of the test cases chosen}
\resizebox{0.48\textwidth}{!}{
\begin{tabular}{lrrrrrrr}
\hline
         & \multicolumn{1}{l}{Load interval} & \multicolumn{2}{l}{VarLower} & \multicolumn{2}{l}{AllSets} & \multicolumn{2}{l}{Best Set} \\
Test case &
  \multicolumn{1}{l}{\begin{tabular}[c]{@{}l@{}} around default \\  value \end{tabular}} &
  \multicolumn{1}{l}{DB} &
  \multicolumn{1}{l}{Sets} &
  \multicolumn{1}{l}{DB} &
  \multicolumn{1}{l}{Sets} &
  \multicolumn{1}{l}{DB} &
  \multicolumn{1}{l}{Sets} \\ \hline
9-bus    & $[50\%,150\%]$                & 13263          & 4           & 13273         & 2           & 13265          & 3           \\
39-bus   & $[75\%,125\%]$                & 13896          & 18          & 13277         & 13          & 13350          & 11          \\
89-bus   & $[90\%,110\%]$                & 17133          & 42          & 18030         & 52          & 13565          & 19          \\
1354-bus & $[75\%,100\%]$                & 19011          & 82          & 21190         & 166          & 14204          & 36          \\
2869-bus & $[75\%,100\%]$                & 29649          & 282         & 29703         & 367         & 18889          & 81         \\ \hline
\\
\end{tabular}}
\label{tab:cases}
\end{table}

\subsection{Baseline Methods}
\subsubsection{KKTs and Big-M}
As mentioned in Section~\ref{sec:Blvl}, the first baseline method used for comparison is obtained by directly solving the MILP in Equations~\eqref{milp_obj}~to~\eqref{milp_bin}. Two versions of this method are implemented. In the first one, later referred as \textit{BigM Long}, the method is left running as long as necessary for it to complete the calculations, within the limit of 15 minutes.
This limit needs to be set because the method is tested for 5'000 scenarios and in the situation where we have many intractable instances, it would not be possible to obtain results in a reasonable amount of time. Moreover, in the case of a strategic generator preparing its bids for the day-ahead market, obtaining an output for each hourly bid within 15 minutes or less would probably be a requirement. In the cases where the runtime limit is reached, the latest feasible solution obtained by the MILP solver is returned, if it is available. The second version of this method will be called \textit{BigM Short}. In this version, the first feasible solution obtained by the MILP solver is returned and the execution stops.

The choice of $M$ is an important factor to take into consideration~\cite{kleinert2019lunch, pineda2019bigm}. Too small, it can interfere with the physics of the model, too big it can lead to numerical ill-conditioning. Here we take advantage of the database generation to set $M$. First, two constants are defined: $M_p$ for the primal constraints and $M_d$ for the dual variables. The choice of $M_p$ is straightforward as it is linked with the existing bounds on the primal variables. However, the choice of $M_d$ is not straightforward as it is complicated to evaluate and bound the dual variables. When solving the DCOPF in the database creation process, the maximum values of the dual variables are retrieved and stored to be later assigned to $M_d$ in the solving process. In order to have some safety margin, we used 10 times these maximum values both for $M_p$ and $M_d$.

\subsubsection{PADM}
A promising method from the literature, the PADM introduced in~\cite{kleinert2021computing}, is also used to compare with our method. As the problem of the strategic generator is not of the form used in this work, we had to adjust the PADM method in order to use it for our problem, as we detail in~\cite{prat2020}. Here again, two versions are implemented. In the first one, \textit{PADM Long}, the method runs until convergence, within the same limit of 15 minutes. In \textit{PADM Short}, the method runs for the same duration as the method \textit{AllSets}, completing the ongoing iteration when the time limit is reached, and the latest values for the variables is kept.

With the PADM too, there are parameters to be determined, which impact the solution efficiency. Those are the initial value of the penalty $\eta_0$ and the penalty increment $\eta^\text{step}$. We observed that if $\eta_0$ is too small, the method can converge to a solution that is not feasible for the original problem. The value of $\eta^\text{step}$ can have a big impact on the number of iterations and if chosen too big, it can impede convergence. The values finally chosen were determined by trial-and-error. These are $\eta_0=1$ and $\eta^\text{step}=0.1$ for 9, 39 and 89-bus systems and $\eta^\text{step}=1$ for 1354 and 2869-bus systems.

\subsection{Modeling Parameters}

\paragraph{Databases generation}
The parameters used in \textit{DiscoverMass} algorithm for the database generation are the ones suggested in~\cite{misra2018learning}. The resulting number of points in the databases is given in Table~\ref{tab:cases}. As the number of active sets varies depending on the method (one with \textit{BestSet}, $\text{n}_\text{A}$ with \textit{AllSets} and $\text{n}_\text{V}$ with \textit{VarLower} in Figure~\ref{fig:flowcharts}), so does the number of points for a given system. The loads applied are selected randomly with equal probability, and independently for each bus, in an interval [$(1-x^\text{m})P^\text{d}_i$,$(1+x^\text{p})P^\text{d}_i$], where $P^\text{d}_i$ is the default load (from the test case data) and $x^\text{m}$ and $x^\text{p}$ are the percentages of sample range under and above the default load, respectively. All the points which are infeasible for the DCOPF are not included the database. The maximum value of $c^\text{S}$, as defined in Section~\ref{sec:Blvl}, is chosen as:
    \begin{equation}
    c^\text{S,max}=10 \times \max\{c_i, i\in I\}
    \end{equation}
where $\max\{c_i, i\in I\}$ is the cost of the most expensive generator in the system. This way we avoid the situation where $c^\text{S}$ is unbounded, while ensuring that the chosen bound will not impact the resulting dispatch.
For methods \textit{AllSets} and \textit{BestSet}, each sample of load was tested with 10 values of $c^\text{S}$, between $c_1$ and $c^\text{S,max}$. Note that the choice of values for $c^\text{S}$ can be crucial because of the influence they have on which constraints are active.
Here, we take a simplified approach and divide the interval of cost in ten. It could be refined by selecting the $c^\text{S}$ values equal to the costs of the other generators, as we observe that those are values for which the active constraints can change. However, a trade-off must be made between the selection of critical values and the size of the database. The generation of the database is something to be improved in future work. Such a study should also focus on the balance of the database and on the representation of each class, both in the training and testing sets.
Here, the scenarios are generated artificially because there is no existing database of operating points. For a real system, it would be possible to learn from observed data.

\begin{figure}[h!]
  \centering
    {\includegraphics[width=0.75\linewidth]{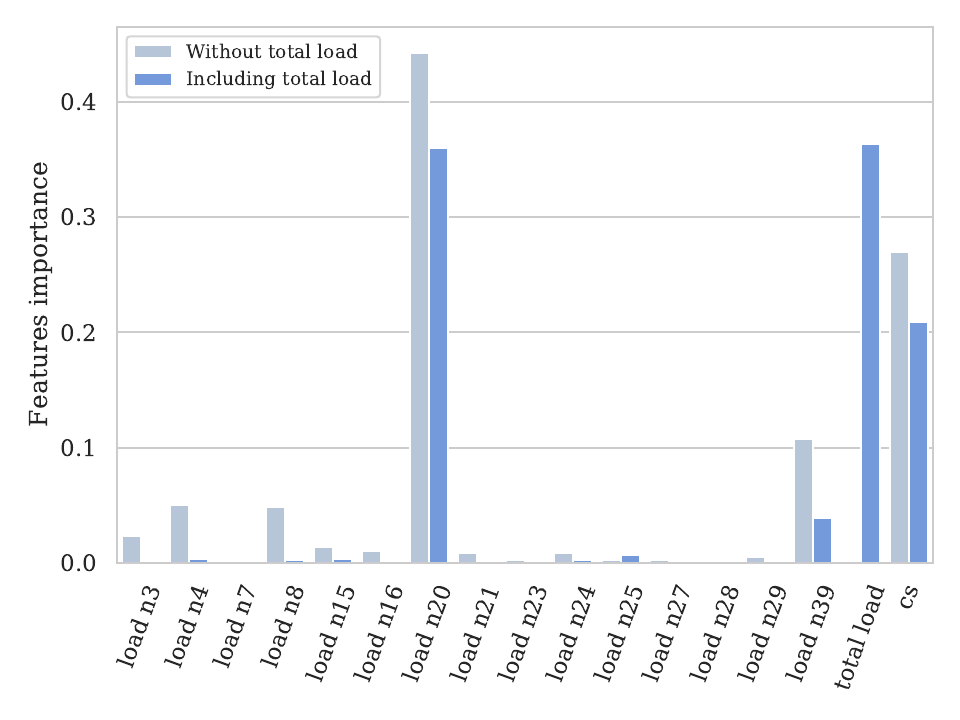}}
  \vspace{-0.5cm}
  \caption{Importance of the different features in the resulting DT. Example with 39-bus system and \textit{VarLower} general DT.}
  \label{fig:features}
\end{figure}

\paragraph{Decision Trees building}
When building a classifier, the features of the data used for the training are of great importance. Here, the main features are the load at each bus and the cost bid by the strategic generator (for Method \textit{VarLower} only). It was found that adding the total load of the system as a feature is very valuable information as well. Figure~\ref{fig:features} illustrates this. The DTs performance improves significantly with this addition as shown in Table~\ref{tab:accuracy}. It has, however, almost no effect on the big 2869-bus system.
To build and evaluate the DTs, the database is randomly separated into training and testing sets, with 70\% of the samples kept for training. The performance is evaluated by the achieved DT accuracy on the test set, which measures the percentage of the test samples that are correctly classified. A randomized search is performed in order to select the hyperparameters that maximize the DT accuracy while avoiding overfitting.

Note that in Table~\ref{tab:accuracy}, we see that even with total load as a feature, the accuracy of the DT can be low for the largest test cases. This will be further explained below in Section~\ref{sec:limits} along with a way to overcome this.
Moreover, we see in this table that the accuracy of the DTs for the 1354-bus system is worse than for the 2869-bus system. To understand why, we must look at the representation of each class in the database. In the case of the method \textit{AllSets}, for the 2869-bus system, the most represented class covers 64\% of the points in the database and the 3 most represented classes cover 90\% of the points in the database. For the 1354-bus system, these numbers drop to 12\% and 32\% respectively. So not only is it harder for the DT to classify the data, but it is also more likely to be wrong for the 1354-bus system than for the 2869-bus system.

\begin{table*}[th] \centering
\caption{Decision tree training time and accuracy depending on the features}
\resizebox{0.8\textwidth}{!}{
\begin{tabular}{lrrrrrrrrr}
 & \multicolumn{3}{c}{Training duration (min)} & \multicolumn{3}{c}{Accuracy without total   load} & \multicolumn{3}{c}{With total load as a   feature} \\
 & VarLower      & AllSets      & BestSet     & VarLower         & AllSets        & BestSet        & VarLower         & AllSets         & BestSet        \\ \hline
9-bus    & 0.35 & 0.01 & 0.01 & 95.0\% & 97.5\% & 97.0\% & 100.0\% & 100.0\% & 100.0\% \\
39-bus   & 3.74 & 0.04 & 0.03 & 86.0\% & 70.3\% & 82.0\% & 96.8\%  & 96.0\%  & 92.6\%  \\
89-bus   & 6.77 & 0.10 & 0.07 & 87.0\% & 72.7\% & 78.1\% & 91.6\%  & 86.4\%  & 91.0\%  \\
1354-bus & 43.17 & 2.48 & 1.56 & 66.0\% & 15.8\% & 21.7\% & 75.2\%  & 25.6\%  & 29.3\%  \\
2869-bus & 277.46 & 16.46 & 8.75 & 52.6\% & 43.5\% & 79.0\% & 52.6\%  & 43.7\%  & 79.2\%  \\ \hline
\end{tabular}}
\\
\label{tab:accuracy}
\end{table*}

\paragraph{Models evaluation}
In order to assess the performance of our methods, 5'000 scenarios of feasible load vectors are randomly generated with a uniform distribution within the same intervals of values that were used for the database generation. The reduced LPs are run in parallel. To ensure feasibility of the solution, we check that the solution of each reduced LP does not violate any constraint of the initial problem. The solution maximizing the profit of the strategic generator is chosen among all feasible solutions.

\subsection{Results}
The different methods are compared in terms of runtime and quality of the solution returned.

\subsubsection{Runtime}
\begin{figure*}[t]
  \centering
    \subfloat[Performance profile for the running times. x-axis represents Eq.~\eqref{ratio}, and is expressed in logarithmic scale.]
    {\includegraphics[width=0.35\linewidth]{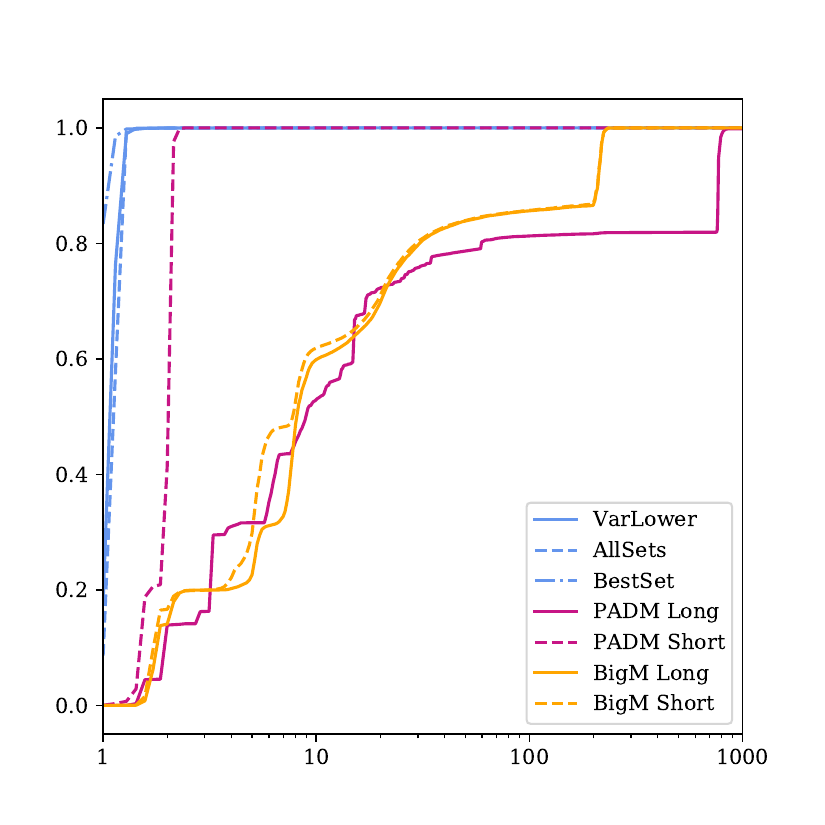}}
    \qquad
    \subfloat[Cumulative distribution function for the relative gaps. x-axis represents Eq.~\eqref{gap}.]
    {\includegraphics[width=0.35\linewidth]{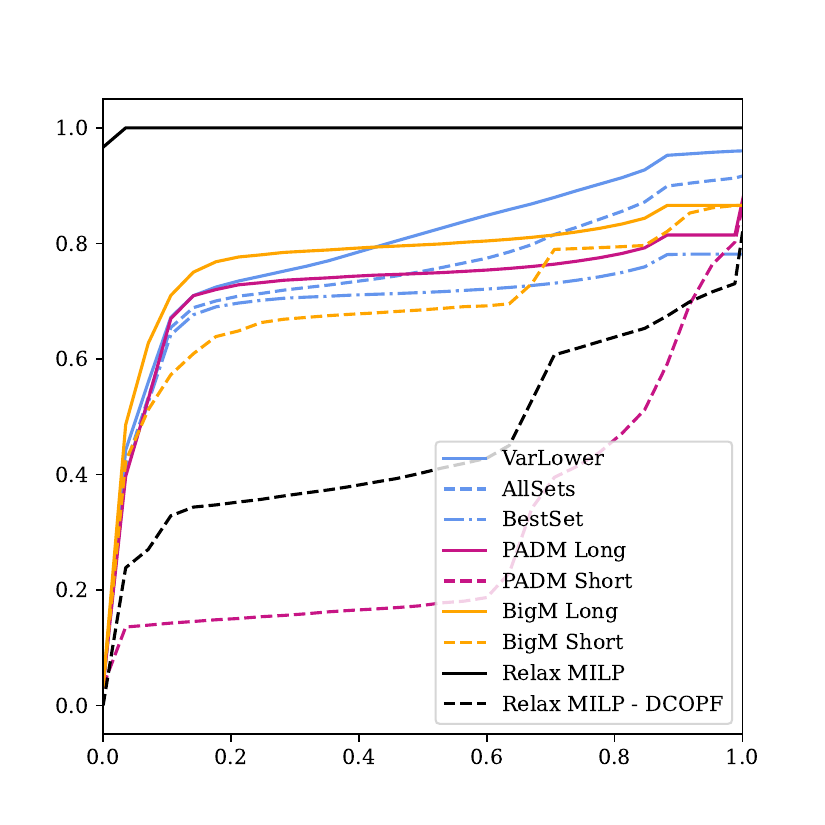}}
  \caption{Comparison of the performance of the different methods in terms of running time and relative gap. y-axis represents the share of instances tested.}
  \label{fig:perf}
\end{figure*}

For the analysis of the temporal component, the performance graph shown in Figure~\ref{fig:perf} on the left, is built as described in~\cite{dolan2002}. This graph shows, for each method, the cumulative distribution function for the performance ratio of method $s$ over all instances $p$. The performance ratio is denoted by $r_{p,s}$, and calculated with:
\begin{equation}
\label{ratio}
r_{p,s} = \frac{t_{p,s}}{\text{min}\{t_{p,s}:s\in S\}} 
\end{equation}
where $S$ represents all the methods that are considered here. In case the method hits the runtime limit, the performance ratio is set to a maximum, $r_M=1000$ in our case.

From this figure, we can see that the methods introduced in this paper manage to return a solution much faster than all the methods that have been used for comparison. More information about the running times per method and per test case is given in Table~\ref{tab:time}.

\begin{table*}[t] \centering
\caption{Results on the running times per method and test case}
\resizebox{0.9\textwidth}{!}{
\begin{tabular}{lrrrrrrrrrr}
 & \multicolumn{2}{c}{9-bus} & \multicolumn{2}{c}{39-bus} & \multicolumn{2}{c}{89-bus} & \multicolumn{2}{c}{1354-bus} & \multicolumn{2}{c}{2869-bus} \\
 & Mean (s)    & Median (s)   & Mean (s)    & Median (s)    & Mean (s)    & Median (s)   & Mean (s)     & Median (s)     & Mean (s)     & Median (s)     \\ \hline
PADM Long  & 0.017 & 0.019 & 5.814 & 0.069 & 3.155 & 0.645 &    809.249    &    900.541    &   52.266    &    26.938   \\
PADM Short & 0.009 & 0.009 & 0.037 & 0.037 & 0.088 & 0.088 & 2.307  & 2.352  &   8.422    &    8.421   \\
BigM Long  & 0.011 & 0.011 & 0.128 & 0.13  & 0.325 & 0.329 & 28.687 & 24.009 &    671.892   &    904.743   \\
BigM Short & 0.011 & 0.011 & 0.097 & 0.096 & 0.3   & 0.34  & 26.996 & 23.112 &    666.819   &    904.624   \\
VarLower    & 0.007 & 0.007 & 0.02  & 0.02  & 0.046 & 0.046 & 1.175  & 1.174  & 4.494 & 4.499 \\
AllSets     & 0.007 & 0.007 & 0.019 & 0.019 & 0.051 & 0.052 & 1.3    & 1.298  & 4.897 & 4.934 \\
BestSet     & 0.006 & 0.006 & 0.018 & 0.018 & 0.043 & 0.043 & 1.186  & 1.184  & 4.258 & 4.259 \\ \hline
\end{tabular}}
\label{tab:time}
\end{table*}

Except for the smaller 9-bus system, for which the big-M method already solves fast, our proposed methods achieve a computation speedup of 6 to 24 times faster than the conventional Big-M method on the tractable instances. For the 2869-bus system, the Big-M method often fails to complete calculation within 15 minutes (which corresponds to a median around 900s) which gives a considerable advantage to our methods for such a large system. Note that for \textit{PADM Long}, the high mean and low median for the running time indicates that some instances do not converge in the selected time limit of 15 minutes.

\subsubsection{Solution quality}
The performance of the different methods in terms of solution quality has been evaluated by comparing the value of the profit for the strategic generator. Instead of directly using the value of the upper-level objective function obtained with the different methods, the profit is calculated with the nodal price returned when solving the DCOPF with the final value of $c^\text{S}$. This is indeed the most accurate way to estimate this profit. The relative gaps $g_{p,s}$ are then calculated with:
\begin{equation}
\label{gap}
g_{p,s} = \frac{q_{p}^*-q_{p,s}}{q_{p}^*} 
\end{equation}
where $q_{p,s}$ is the profit of the strategic generator for instance $p$, with method $s$, and $q_{p}^*$ is the lower bound on the MILP optimal objective function, which is obtained by relaxing the binary constraints in Problem~\eqref{eq:milp}. We check the value of the now integer but originally binary variables. If they are all close enough to 0 or 1, we can conclude that the optimal solution of the relaxed MILP is the solution of the MILP, which is never the case here. In other words, the solution of the relaxed MILP is never in the feasible space of the MILP here. For this reason, we include one more result as an upper bound on the optimal value of the objective function of the MILP, which is obtained by solving the DCOPF with the optimal value of $c^\text{S}$ in the relaxed MILP. Those two supplementary comparisons are included with the names \textit{Relax MILP} and \textit{Relax MILP-DCOPF} respectively.
The resulting cumulative distribution function is shown on the right part of Figure~\ref{fig:perf}. It can be observed that the performance of the three methods we propose in this paper, in terms of solution value, is close to the lower-bound obtained with \textit{Relax MILP}.
In particular, the methods perform much better than \textit{BigM Short} and \textit{PADM Short}. Details on the relative gaps per method and per test case are given in Table~\ref{tab:gap}. This confirms that, when infeasible instances are disregarded, the solution quality is good with our methods.

\begin{table*}[t] \centering
\caption{Results on the relative gaps, comparing to the lower bound and excluding infeasible instances}
\resizebox{0.8\textwidth}{!}{
\begin{tabular}{lrrrrrrrrrr}
 & \multicolumn{2}{c}{9-bus} & \multicolumn{2}{c}{39-bus} & \multicolumn{2}{c}{89-bus} & \multicolumn{2}{c}{1354-bus} & \multicolumn{2}{c}{2869-bus} \\
 & Mean & Median & Mean & Median & Mean & Median & Mean & Median & Mean & Median \\ \hline
PADM Long  & 11.0\%  & 0.0\%  & 28.7\%  & 3.2\% & 5.7\%  & 3.7\%  & 89.9\% & 100.0\% & 8.4\%  & 8.4\%  \\
PADM Short & 80.0\% & 83.1\% & 38.1\% & 4.6\% & 65.4\% & 65.6\% & 101.2\% & 100.0\% & 88.1\%  & 89.3\%  \\
BigM Long  & 11.0\%  & 0.0\%  & 28.7\%  & 3.2\% & 5.7\%  & 3.7\%  & 5.6\% & 3.5\% & 9.3\%  & 7.4\%  \\
BigM Short & 11.0\%  & 0.0\%  & 36.2\% & 11.9\% & 35.9\% & 32.8\%  & 8.7\% & 4.1\% & 13.8\%  & 7.6\%  \\
VarLower    & 11.0\%  & 0.0\%  & 29.1\%  & 4.3\% & 7.3\%  & 3.7\%  & 29.5\% & 14.9\% & 18.6\%  & 9.3\%  \\
AllSets     & 11.0\%  & 0.0\%  & 30.7\%  & 6.3\% & 8.8\%  & 3.8\%  & 43.9\% & 48.9\% & 11.4\%  & 8.4\% \\
BestSet     & 11.0\%  & 0.0\%  & 30.1\%  & 4.9\% & 5.2\%  & 3.7\%  & 5.7\% & 2.9\% & 8.5\% & 8.4\% \\ \hline
\end{tabular}}
\label{tab:gap}
\end{table*}

\subsubsection{Duality Gap of the Lower-Level Problem}
For given values of the strategic cost $c^\text{S}$ and of the lower-level primal and dual optimization variables, the duality gap is calculated as the difference between the primal and dual objective functions. The value of the dual objective function is given by: $ \sum_{i} (\alpha_i P_i^\text{d} + \phi_i^\text{min} P_i^\text{g,min} - \phi_i^\text{max} P_i^\text{g,max}) - \sum_{l} f_{l}^\text{max} (\rho_l^\text{min}+\rho_l^\text{max}) $.

The three methods introduced here are designed in a way that if a feasible solution is found, the duality gap of the lower-level problem will always be zero. In other words, if these methods return a feasible point, it is the optimal solution. The same applies to the big-M methods, and to the PADM if it converges. The \textit{Relax MILP}, on the other hand, can have a duality gap. For example, for the 1354-bus system, the relative duality gap is on average $4.0\%$ with \textit{Relax MILP} and $3.6\%$ with \textit{PADM Short}.

\subsubsection{Limits of the Methods and Improvements} \label{sec:limits}
In Table~\ref{tab:in}, the number of instances that reach the time limit and the number of instances that do not return a feasible solution are given. This confirms that the methods introduced in this paper have the advantage of avoiding intractability. However, the reformulation introduced might be infeasible and unable to return a value for the bidding cost of the strategic generator. The PADM always yields a value, as, even in the case it is not converging, an estimation of the strategic cost is made at each iteration. For the 1354-bus system in particular, \textit{PADM Long} fails to converge for most instances but still gives a value for $c^\text{S}$, which appears to be far from the optimal value as show in Table~\ref{tab:gap}.

\begin{table*}[t] \centering
\caption{Results on intractability and infeasibility}
\resizebox{1\textwidth}{!}{
\begin{tabular}{lrrrrrrrrrr}
        & \multicolumn{2}{c}{9-bus} & \multicolumn{2}{c}{39-bus} & \multicolumn{2}{c}{89-bus} & \multicolumn{2}{c}{1354-bus} & \multicolumn{2}{c}{2869-bus} \\
        & Time limit   & Infeasible  & Time limit   & Infeasible   & Time limit  & Infeasible   & Time limit   & Infeasible     & Time limit   & Infeasible     \\ \hline
PADM Long           & 0 & 0 & 31 (0.6\%) & 0 & 13 (0.3\%) & 0           & 4472 (89\%)  & 0            & 44 (0.9\%)            & 0             \\
PADM Short          & 0 & 0 & 0          & 0 & 0  & 0           & 0 & 0            & 0             & 0             \\
BigM Long           & 0 & 0 & 0          & 0 & 0  & 0           & 0 & 0            & 3186 (63.7\%) & 3186 (63.7\%) \\
BigM Short          & 0 & 0 & 0          & 0 & 0  & 0           & 0 & 0            & 0             & 0             \\
VarLower             & 0 & 0 & 0          & 0 & 0  & 48 (1.0\%)  & 0 & 318 (6.4\%)  & 0             & 442 (8.8\%)             \\
AllSets              & 0 & 0 & 0          & 0 & 0  & 132 (2.6\%) & 0 & 682 (13.6\%) & 0             & 926 (18.5\%)  \\
BestSet & 0            & 0           & 0            & 239 (4.8\%)  & 0           & 418 (8.4\%)  & 0            & 3544 (70.9\%)  & 0            & 1089 (21.8\%)  \\
VarLower (augmented) & - & - & -          & - & -  & -           & 0 & 22 (0.4\%)            & 0             & 284 (5.7\%)             \\
AllSets (augmented)  & - & - & -          & - & -  & -           & 0 & 5 (0.1\%)           & 0             & 63 (1.3\%)             \\
BestSet (augmented)  & - & - & -          & - & -  & -           & 0 & 48 (1.0\%)            & 0             & 13 (0.3\%)             \\ \hline
\end{tabular}}
\label{tab:in}
\end{table*}

\begin{figure*}[t]
  \centering
    \subfloat[1354 bus]
    {\includegraphics[width=0.47\linewidth]{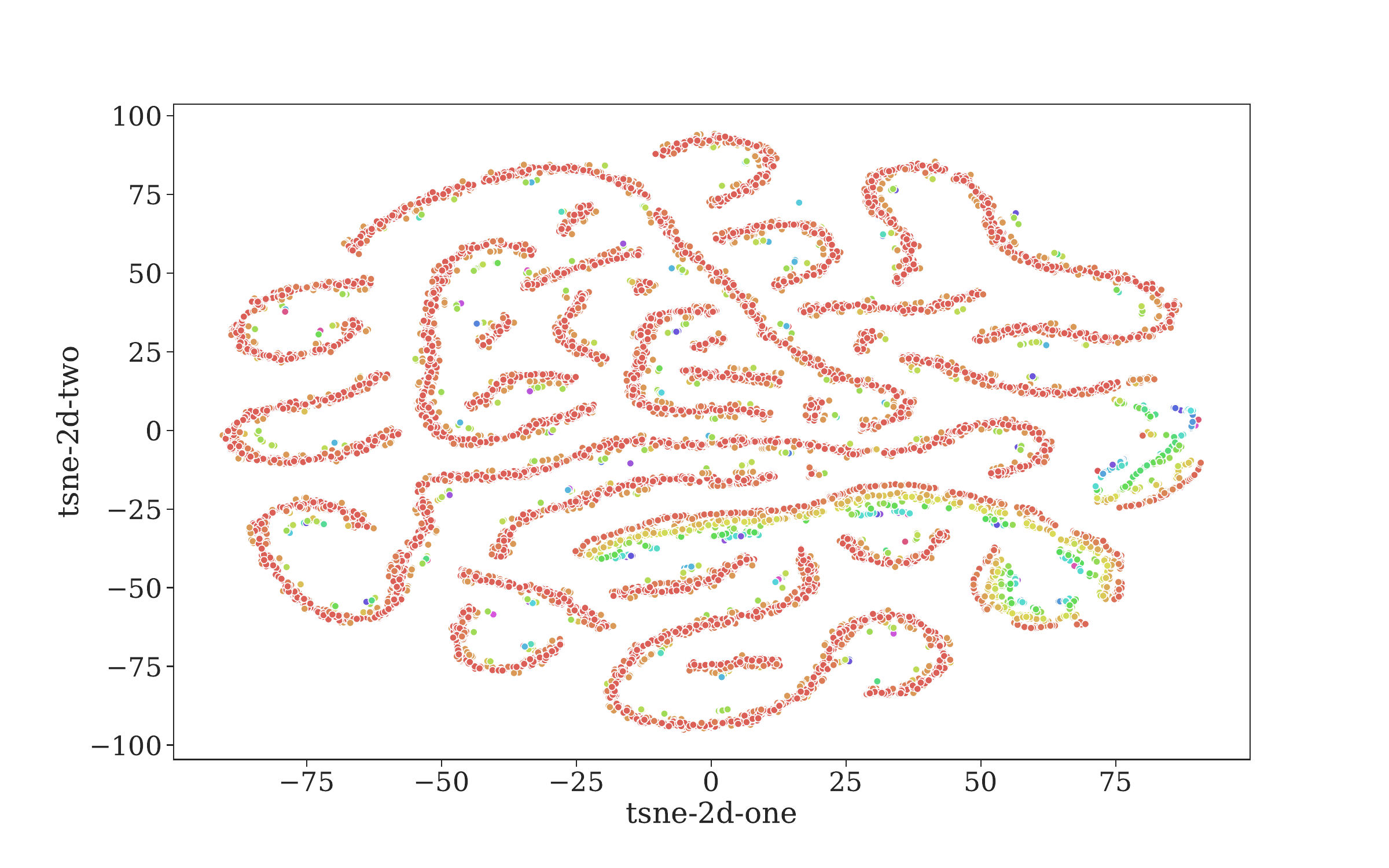}}
    \qquad
    \subfloat[2869 bus]
    {\includegraphics[width=0.47\linewidth]{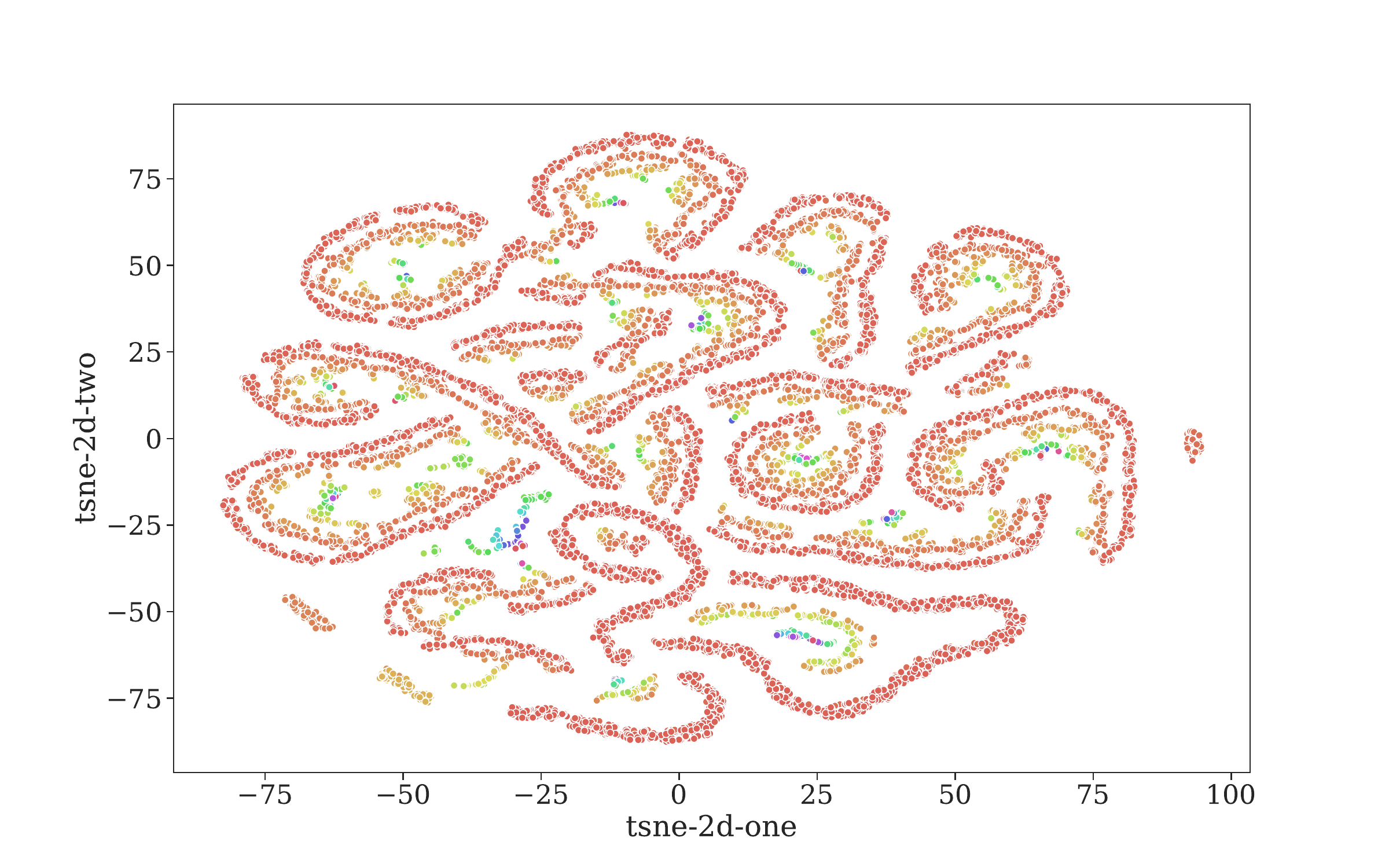}}
  \caption{t-SNE plots for the complex cases, on the databases generated for \textit{VarLower}. Each color corresponds to a possible set of active constraints.}
  \label{fig:tsne}
\end{figure*}

For the complex systems studied (1354-bus and 2869-bus) the incidence of infeasible occurrences is high. This derives from the fact that the DTs have a low accuracy. Figure~\ref{fig:tsne} presents t-SNE plots for the databases with the method \textit{VarLower}. The t-SNE technique allows to represent a high-dimensional set of data in 2 dimensions. It shows points that are related by clustering them. Ideally, there is a cluster for each class of the data represented. Further information on t-SNE can be found in~\cite{vandermaaten2008tsne}. Here, each class corresponds to a set of active constraints. However, we cannot identify proper clusters, which means that the classes of our problem are very similar. It is then difficult for the classifier to distinguish how to best separate the data, while avoiding over-fitting. But similar classes will be classified in the same area of the created decision tree. This idea is applied to extract extra information from the DT and significantly increase the accuracy for these systems. When applying the DT to a given load, the corresponding leaf is extracted and its parent is identified. Then, the classes of all training data samples that would be classified to the parent node of this leaf are retrieved. As a result, more LPs are formulated but a good performance is recovered. The three last rows of Table~\ref{tab:in} indicate that the number of infeasible instances decreases significantly when this technique is applied.

For the 2869-bus system, from 64\% of intractable cases with the Big-M method, we manage to reduce them to about 20\% of infeasible cases with our methods applied directly and as low as 0.3\% when we extract more information from the learning process of the DT with our augmented method. 

\subsubsection{Comparison of the three Methods}
We have seen that the three methods introduced have similar performance but there are some differences that we will develop here. With the method \textit{BestSet}, less active sets are identified, as shown in Table~\ref{tab:cases}. As a consequence, the database is smaller (and thus faster to build) and it is faster to train the DT, as we see in Table~\ref{tab:accuracy}. This method gives the best results in terms of distance to optimality. However, it results more often to an infeasible model, which is coherent with the fact that only one LP is solved in that case.
With the method \textit{VarLower}, more DTs need to be trained, which takes more time.
This method results less often than the others in infeasible instances. However, it cannot be applied to bilevel problems where there is more than one decision variable passed on to the lower-level problem.
Finally, with the method \textit{AllSets}, a lot of sets of sets of active constraints can be identified and it is complicated to ensure that they are all properly represented in the database. However, since multiple LPs are run, it gives better results than \textit{BestSet} in terms of infeasibility.

The choice of one method over another thus depends on one's focus. If more interested in retrieving the optimal solution, \textit{BestSet} should be prefered. On the other hand, if the focus is to find a feasible solution, \textit{VarLower} will perform better. Alternatively, \textit{AllSets} can be used, especially for other bilevel problems in which more than one decision variable are passed on to the lower-level problem.

\section{Conclusion}
\label{sec:ccl}
This paper uses machine learning to introduce efficient approaches that solve linear bilevel problems, boosting their runtime and solution quality. We propose three methods, some of them highly parallelizable, that use decision trees to learn the active sets of the lower-level problem and apply them to the problem of a strategic generator optimizing its bids for the electricity market. Linear bilevel problems are most often converted to Mixed Integer Linear Programs (MILP) by replacing the lower-level problem with its Karush--Kuhn--Tucker conditions (KKTs), in order to solve them. Our goal in this paper is to completely eliminate the use of the binary variables by learning the active sets and solve a single or a small number of Linear Programs (LPs) instead.

Contrary to existing machine learning methods for optimization, to the best of our knowledge, this is the first paper that considers optimization problems with two distinct characteristics: besides effectively treating MILPs -- and not LPs as most machine learning methods applied on optimization problems so far -- the main challenge is that the decision variables of the upper-level problem shall not intervene with the active set classification process of the lower-level problem, although the upper-level decision variables are indeed parameters of the lower-level problem.

We apply our methods to systems with up to 2'869 buses, and we compare them with the most promising existing approaches to solve bilevel problems: the big M reformulation, which is the most common approach to treat the non-linear complementarity constraints, and the Penalty Alternating Direction Method (PADM), recently introduced in~\cite{kleinert2021computing}. Our methods are shown to be 6 to 24 times faster, while they maintain very good solution quality, comparable with existing methods. 
More importantly, although we observe a trade-off between the decision tree prediction accuracy and the computational complexity of the problem (the number of LPs solved), our methods have been shown to retrieve good solutions to problems where existing methods, such as the PADM, fail to converge. At the same time, they also appear to perform well in large test cases. As a matter of fact, we obtain the most encouraging results for very large systems, such as the 2869-bus system, where existing methods fail to return a solution within 15 minutes for 64\% of the test cases. Finally, an additional considerable advantage of the methods proposed in this paper is that they do not introduce any parameters in the online stage, that would have to be appropriately tuned for a good performance.

As far as power system problems are concerned, considering that we already use training data that include a wide range of possible load and generation realizations, the methods proposed in this paper can easily accommodate uncertainty from various sources, and especially renewable generation, while they can also be extended to a bilevel program with a stochastic lower-level problem. Future work should focus on providing guarantees for the accuracy of the Decision Trees, when it comes to the classification of the active sets; this shall enhance the performance of the proposed methods. Research is also required to investigate the desirable properties of the database and how to achieve them, and to extend the application of such approaches to more complex bilevel problems.

\bibliographystyle{IEEEtran}

\bibliography{Bib}

\begin{IEEEbiography}[{\includegraphics[width=1in,height=1.25in,clip,keepaspectratio]{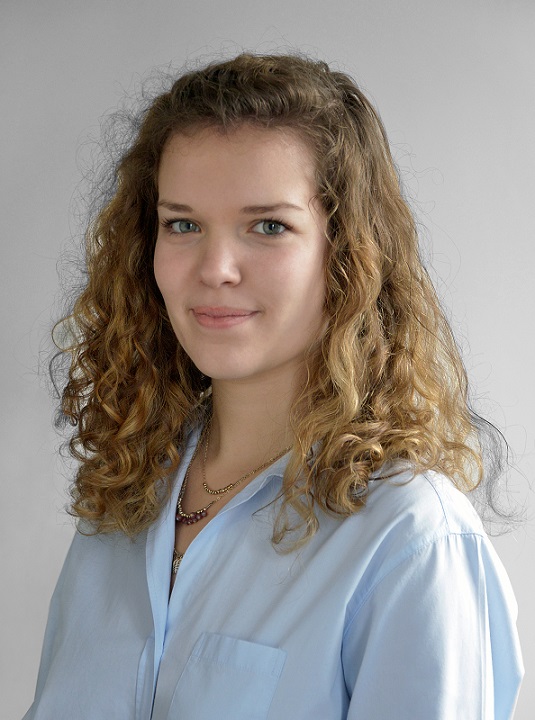}}]{Eléa Prat}
received her M.S. in sustainable energy from the Technical University of Denmark (DTU) in 2017. She worked as a research assistant in the electrical engineering department of DTU (2019-2021). She is now doing her Ph.D. in the management department of DTU, in operations research. Her research interests include optimization, electricity markets, multi-scale decision making, bilevel programming.

\end{IEEEbiography}

\begin{IEEEbiography}[{\includegraphics[width=1in,height=1.25in,clip,keepaspectratio]{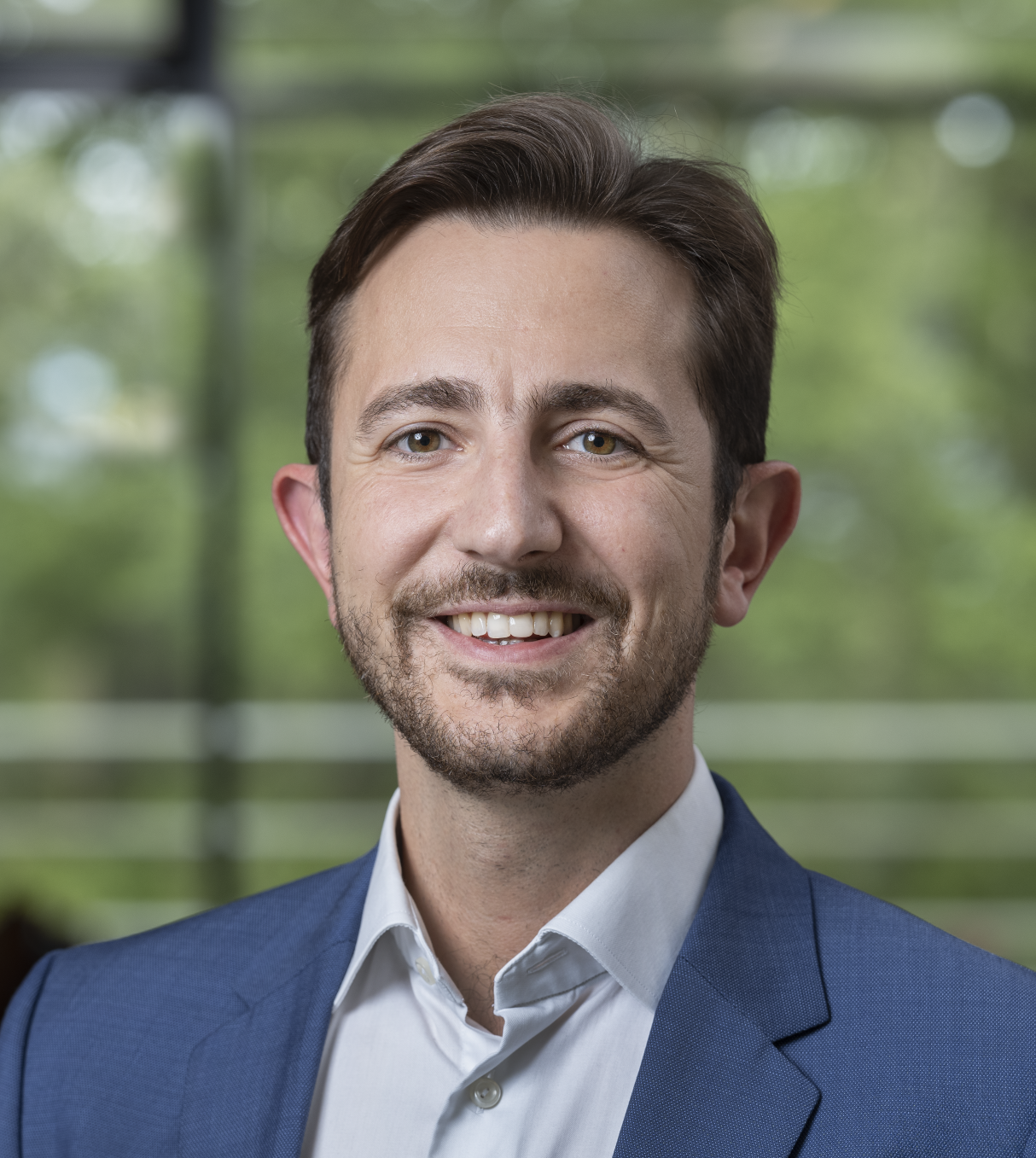}}]{Spyros Chatzivasileiadis}
(S’04, M’14, SM’18) is the Head of Section Power Systems and an Associate Professor at the Technical University of Denmark (DTU). Before that he was a postdoctoral researcher at the Massachusetts Institute of Technology (MIT), USA and at Lawrence Berkeley National Laboratory, USA. Spyros holds a PhD from ETH Zurich, Switzerland (2013) and a Diploma in Electrical and Computer Engineering from the National Technical University of Athens (NTUA), Greece (2007). He is currently working on machine learning applications for power systems, and on power system optimization, dynamics, and control of AC and HVDC grids. Spyros is the recipient of an ERC Starting Grant in 2020.

\end{IEEEbiography}

\appendix
\subsection{PADM Formulation for the Strategic Generator Problem}
\label{ap:PADM}

\subsubsection{Dual of Lower-Level Problem}
The dual of the lower-level problem given in Equations~\eqref{strat_lobj} to \eqref{strat_ref} is:

\begin{footnotesize}
\begin{subequations}
\begin{alignat}{2}
\label{eq:dual_obj} \underset{\alpha, \phi, \rho, \gamma}{\text{min}} \quad & \sum_{i} (\phi_i^\text{max} P_{i}^\text{g,max} - \phi_i^\text{min} P_{i}^\text{g,min}) + \sum_{l} f_{l}^\text{max} (\rho_l^\text{min}+\rho_l^\text{max})\\
\label{eq:dual_Pg1} \quad \text{s.t.} \quad & c^\text{S} - \alpha_{i=1} - \phi_{i=1}^\text{min} + \phi_1^\text{max} = 0 \\
& \label{eq:dual_Pg} c_i - \alpha_i - \phi_i^\text{min} + \phi_i^\text{max} = 0 , \quad \forall i \neq 1\\
\nonumber & \sum_{l, i=\text{from}_l} B_l (\alpha_i - \alpha_{i=\text{to}_l} - \rho_l^\text{min} + \rho_l^\text{max}) \\
\label{eq:dual_theta} & + \sum_{l, i = \text{to}_l} B_l (\alpha_i - \alpha_{i=\text{from}_l} + \rho_l^\text{min} - \rho_l^\text{max}) = 0, \quad \forall i \neq \text{ref}\\
\nonumber & \sum_{l, \text{ref} = \text{from}_l} B_l (\alpha_{i=\text{ref}} - \alpha_{i=\text{to}_l} - \rho_l^\text{min} + \rho_l^\text{max}) \\
\label{eq:dual_ref} & + \sum_{l, \text{ref} = \text{to}_l} B_l (\alpha_{i=\text{ref}} - \alpha_{i=\text{from}_l} + \rho_l^\text{min} - \rho_l^\text{max}) + \gamma = 0\\
\label{eq:dual_posi} & \phi_i^\text{min} , \phi_i^\text{max} \geq 0, \quad \forall i\\
\label{eq:dual_posl} & \rho_l^\text{min} , \rho_l^\text{max} \geq 0, \quad \forall l\\
\end{alignat}
\end{subequations}
\end{footnotesize}

\subsubsection{Reformulation as one-level problem}
Using the dual of the lower-level problem and the strong-duality theorem, the bilevel problem in Equations~\eqref{strat_uobj}
to \eqref{strat_ll} can be reformulated as a single-level problem:
\begin{footnotesize}
\begin{subequations}
\begin{alignat}{2}
\label{eq:one_obj} \underset{c^\text{S}, P^\text{g},\theta,\alpha, \phi, \rho, \gamma}{\text{min}} \quad & c_{i=1} P_{i=1}^\text{g}  - \alpha_{i=1} P_{i=1}^\text{g}  \\
\label{eq:one_cap} \quad \text{s.t.} \quad & c_{i=1} \leq  c^\text{S}  \leq c^\text{S,max} \\
\label{eq:one_bal} & P_i^\text{g}-P_i^\text{d}- \sum_{l, i\in l} B_l \Delta\theta_{l} = 0, \quad  \forall i \\
& \label{eq:one_gen} P_i^\text{g,min} \leq P_i^\text{g} \leq P_i^\text{g,max} , \quad \forall i\\
& \label{eq:one_line} -f_{l}^\text{max} \leq  B_l \Delta\theta_{l} \leq f_{l}^\text{max} , \quad \forall l\\
\label{eq:one_setref} & \theta_\text{ref}=0 \\
& \label{eq:one_Pg1} c^\text{S} - \alpha_{i=1} - \phi_{i=1}^\text{min} + \phi_{i=1}^\text{max} = 0 \\
& \label{eq:one_Pg} c_i - \alpha_i - \phi_i^\text{min} + \phi_i^\text{max} = 0 , \quad \forall i \neq 1\\
\nonumber & \sum_{l, i = \text{from}_l} B_l (\alpha_i - \alpha_{i = \text{to}_l} - \rho_l^\text{min} + \rho_l^\text{max}) \\
\label{eq:one_theta} & + \sum_{l, i = \text{to}_l} B_l (\alpha_i - \alpha_{i = \text{from}_l} + \rho_l^\text{min} - \rho_l^\text{max}) = 0, \quad \forall i \neq \text{ref}\\
\nonumber & \sum_{l, \text{ref} = \text{from}_l} B_l (\alpha_{i=\text{ref}} - \alpha_{i =\text{to}_l} - \rho_l^\text{min} + \rho_l^\text{max}) \\
\label{eq:one_ref} & + \sum_{l, \text{ref} = \text{to}_l} B_l (\alpha_{i=\text{ref}} - \alpha_{i =\text{from}_l} + \rho_l^\text{min} - \rho_l^\text{max}) + \gamma = 0\\
\label{eq:one_posi} & \phi_i^\text{min} , \phi_i^\text{max} \geq 0, \quad \forall i\\
\label{eq:one_posl} & \rho_l^\text{min} , \rho_l^\text{max} \geq 0, \quad \forall l\\
\nonumber & - c^\text{S} P_{i=1}^\text{g} - \sum_{i\neq 1} c_i P_{i}^\text{g} + \sum_{i} (P_{i}^\text{d} \alpha_i + P_{i}^\text{g,min} \phi_i^\text{min}\\
\label{eq:one_duality} &  -  P_{i}^\text{g,max} \phi_i^\text{max}) - \sum_{l} f_{l}^\text{max} (\rho_l^\text{min}+\rho_l^\text{max}) \geq 0
\end{alignat}
\end{subequations}
\end{footnotesize}

\subsubsection{Penalty Problem}
In the previous formulation, the complicating constraint~\eqref{eq:one_duality} is relaxed, using the penalty $\eta$:
\begin{footnotesize}
\begin{subequations}
\begin{alignat}{2}
\nonumber \underset{c^\text{S}, P^\text{g},\theta,\alpha, \phi, \rho, \gamma}{\text{min}} \quad & c_{i=1} P_{i=1}^\text{g}  - \alpha_{i=1} P_{i=1}^\text{g}  + \eta [c^\text{S} P_{i=1}^\text{g} \\
\nonumber & + \sum_{i\neq 1} c_i P_{i}^\text{g} - \sum_{i} (P_{i}^\text{d} \alpha_i + P_{i}^\text{g,min} \phi_i^\text{min}  \\
\label{eq:penalty_obj}  & -  P_{i}^\text{g,max} \phi_i^\text{max}) + \sum_{l} f_{l}^\text{max} (\rho_l^\text{min}+\rho_l^\text{max})]\\
\label{eq:penalty_all} \quad \text{s.t.} \quad & \eqref{eq:one_cap} -  \eqref{eq:one_posl}
\end{alignat}
\end{subequations}
\end{footnotesize}

\subsubsection{Subproblems}
From the penalty problem, subproblems are formulated and solved alternatively.
The first is obtained by fixing the primal variables of lower-level problem:
\begin{footnotesize}
\begin{subequations}
\begin{alignat}{2}
\nonumber \underset{c^\text{S},\alpha, \phi, \rho, \gamma}{\text{min}} \quad & - \alpha_1 \mathbf{\bar{P_{i=1}^\text{g}}}  + \eta [ c^\text{S} \mathbf{\bar{P_{i=1}^\text{g}}} - \sum_{i} (P_{i}^\text{d} \alpha_i + P_{i}^\text{g,min} \phi_i^\text{min}\\
\label{eq:sub1_obj}  &  -  P_{i}^\text{g,max} \phi_i^\text{max}) + \sum_{l} f_{l}^\text{max} (\rho_l^\text{min}+\rho_l^\text{max})]\\
\label{eq:sub1_cap} \quad \text{s.t.} \quad & c_1 \leq  c^\text{S}  \leq c^\text{S,max} \\
& \label{eq:sub1_Pg1} c^\text{S} - \alpha_1 - \phi_1^\text{min} + \phi_1^\text{max} = 0 \\
& \label{eq:sub1_Pg} c_i - \alpha_i - \phi_i^\text{min} + \phi_i^\text{max} = 0 , \quad \forall i \neq 1\\
\nonumber & \sum_{l, i =\text{from}_l} B_l (\alpha_i - \alpha_{i =\text{to}_l} - \rho_l^\text{min} + \rho_l^\text{max}) \\
\label{eq:sub1_theta} & + \sum_{l, i =\text{to}_l} B_l (\alpha_i - \alpha_{i =\text{from}_l} + \rho_l^\text{min} - \rho_l^\text{max}) = 0, \quad \forall i \neq \text{ref}\\
\nonumber & \sum_{l, \text{ref} =\text{from}_l} B_l (\alpha_\text{ref} - \alpha_{i =\text{to}_l} - \rho_l^\text{min} + \rho_l^\text{max}) \\
\label{eq:sub1_ref} & + \sum_{l, \text{ref} =\text{to}_l} B_l (\alpha_\text{ref} - \alpha_{i =\text{from}_l} + \rho_l^\text{min} - \rho_l^\text{max}) + \gamma = 0\\
\label{eq:sub1_posi} & \phi_i^\text{min} , \phi_i^\text{max} \geq 0, \quad \forall i\\
\label{eq:sub1_posl} & \rho_l^\text{min} , \rho_l^\text{max} \geq 0, \quad \forall l
\end{alignat}
\end{subequations}
\end{footnotesize}
The second subproblem is obtained by fixing the variables of upper-level and the dual variables of lower-level:
\begin{footnotesize}
\begin{subequations}
\begin{alignat}{2}
\label{eq:sub2_obj} \underset{P^\text{g},\theta}{\text{min}} \quad & c_{i=1} P_{i=1}^\text{g}  - \mathbf{\bar{\alpha_1}} P_{i=1}^\text{g}  + \eta [\mathbf{\bar{c^\text{S}}} P_{i=1}^\text{g} + \sum_{i\neq 1} c_i P_i^\text{g}]\\
\label{eq:sub2_bal} & P_{i}^\text{g}-P_{i}^\text{d}- \sum_{l, i\in l} B_l \Delta\theta_{l} = 0, \quad  \forall i \\
& \label{eq:sub2_gen} P_i^\text{g,min} \leq P_i^\text{g} \leq P_i^\text{g,max} , \quad \forall i\\
& \label{eq:sub2_line} -f_{l}^\text{max} \leq  B_l \Delta\theta_{l} \leq f_{l}^\text{max} , \quad \forall l\\
\label{eq:sub2_setref} & \theta_\text{ref}=0
\end{alignat}
\end{subequations}
\end{footnotesize}

\end{document}